\pdfoutput=1
\documentclass[11pt,a4paper]{article}
\usepackage[T1,OT1]{fontenc}
\usepackage{amsmath,amssymb,amsthm,booktabs,geometry,graphicx,hyperref}
\hypersetup{hidelinks,pdfauthor={Robert Małysz},pdftitle={Dimensions of graphs of vector-valued affine fractal interpolation functions and of their linear images},pdfkeywords={fractal interpolation functions, self-affine sets, Hausdorff dimension, box dimension, projection entropy},pdfsubject={Fractal graph and range dimensions, projections and Bernoulli measures}}
\newtheorem{theorem}{Theorem}
\newtheorem{corollary}{Corollary}
\newtheorem{lemma}{Lemma}[section]
\newtheorem{proposition}[lemma]{Proposition}
\newtheorem{littheorem}{Theorem}

\theoremstyle{definition}
\newtheorem{definition}[lemma]{Definition}
\newtheorem{remark}[lemma]{Remark}
\newtheorem{example}[lemma]{Example}
\theoremstyle{plain}
\newcommand{\R}{\mathbb R}
\newcommand{\E}{\mathbb E}
\newcommand{\B}{\mathrm B}
\newcommand{\HH}{\dim_{\mathrm H}}
\newcommand{\hw}{\mathsf h}
\title{Dimensions of graphs of vector-valued affine fractal interpolation functions and of their linear images}
\author{{\fontencoding{T1}\fontfamily{lmr}\selectfont Robert Małysz}\\[0.35em]
{\small Faculty of Mathematics and Information Science, Warsaw University of Technology}\\
{\small Koszykowa 75, 00-662 Warsaw, Poland}\\
{\small\href{mailto:robert.malysz@pw.edu.pl}{\texttt{robert.malysz@pw.edu.pl}}\qquad ORCID: \href{https://orcid.org/0009-0008-2687-7632}{\nolinkurl{https://orcid.org/0009-0008-2687-7632}}}}\date{}
\begin{document}\maketitle
\begin{abstract}
We determine the common Hausdorff and box dimension of every linear image of the graph of a vector-valued affine fractal interpolation function with positive vertical scaling factors. This includes the graph itself, the range of the function, and the graphs of its coordinate functions. The formula extends to signed factors under a rank condition on the branches with nonzero factors. If all factors are nonzero and their sum differs from one, it holds for every choice of interpolation points and every linear image. The dimension depends only on the horizontal and vertical contraction ratios, the rank of the projected deviations from the endpoint line, and the affine dimension of the projected interpolation points. The proof uses entropy estimates for tilted projections of homogeneous triangular systems, uniform over bounded sets of slopes. For these systems, we show that self-affine measures with positive weights and their tilted projections are exact dimensional, and we compute their dimensions. We give examples in which the rank condition fails and the Hausdorff and box dimensions differ. We also give counterexamples to Massopust's box-dimension formula for systems with diagonal vertical parts.
\end{abstract}

\section{Introduction}
For a vector-valued affine fractal interpolation function (FIF) $f:[0,1]\to\R^n$, the number of coordinates can exceed the number of independent directions in which $f$ departs from its endpoint line. Subtract that line from the interpolation values and let $r$ be the dimension of the span of the residuals. Write $b_i$ for the horizontal widths and $s_i$ for the vertical scaling factors. When all $s_i$ are positive, these parameters and the rank determine the graph dimension. Theorem~\ref{main} determines the dimension of every linear image, including maps that mix the time and value coordinates. For each image, only the branch contractions and two ranks computed from the projected interpolation points enter the formula.

For scalar affine interpolation, B\'ar\'any, Rams and Simon~\cite[Theorem~2.1]{brs} obtain the graph dimension $d$ from $\sum_i |s_i|b_i^{d-1}=1$ when the vertical factors are nonzero, $\sum_i|s_i|>1$ and their nondegeneracy condition on the branch shears holds. Their proof uses the planar result of B\'ar\'any, Hochman and Rapaport~\cite{bhr}. Theorem~\ref{main} extends this graph formula to arbitrary residual rank $r$. When $\sum_i|s_i|>1$, let $q>1$ be the unique solution of $\sum_i|s_i|^q=1$. A new regime appears when $q\le r$. Under the hypotheses for equality, the graph dimension is then $q$ and, at fixed rank, is independent of the horizontal widths $b_i$. This regime first occurs at $r=2$.

\begin{figure}[htbp]
\centering
\includegraphics[width=.96\linewidth]{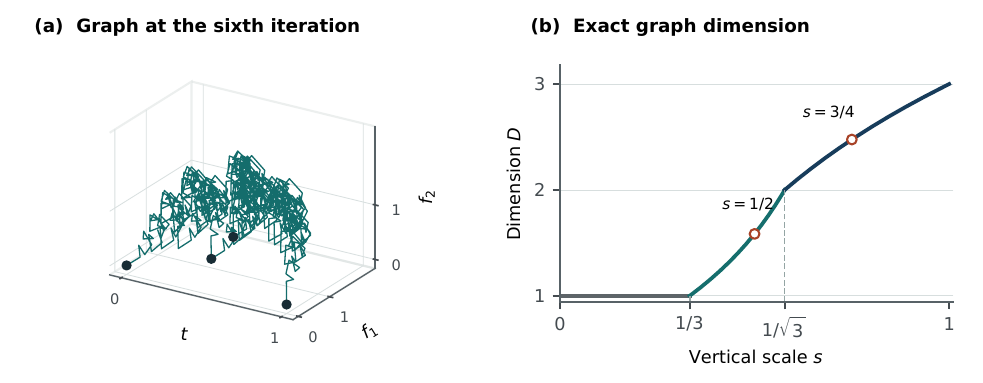}
\caption{Three equal branches with residual rank two. Left: the sixth polygonal iterate at $s=3/4$, started from the zero endpoint interpolant. Dots mark the four interpolation points on the graph. Right: the dimension given by Theorem~\ref{main} for $0<s<1$. It equals one through $s=1/3$, reaches two at $s=1/\sqrt3$, and tends to three as $s$ tends to one. Open circles mark the two computed examples.}
\label{fig:threebranch}
\end{figure}

Take three branches of equal width. The interpolation values are $(0,0),(1,0),(0,1),(0,0)$ at $0,1/3,2/3,1$. The endpoint line is zero and $r=2$. For a common vertical factor $s=1/2$ the root is $q=\log_2 3<2$ and the graph has this dimension. At $s=3/4$ the same interpolation data give three distinct dimensions. The vector graph has dimension $3+2\log_3(3/4)\approx2.476$, its range in the plane has dimension $2$, and the graph of either scalar coordinate has dimension $2+\log_3(3/4)\approx1.738$. In all three cases the Hausdorff and Minkowski dimensions coincide. Figure~\ref{fig:threebranch} shows the vector graph and the change of regime as $s$ varies.

After subtracting the endpoint line, the linear part of branch $i$ sends $(t,z)$ to $(b_it,t\Delta_i+s_i z)$ on $\R\times E$, where $E$ is the residual space and $\Delta_i$ is the increment between consecutive residual interpolation values. This increment is the shear vector of the normalized branch. Every vertical subspace $\{0\}\times V$, $V\le E$, is invariant. For $r>0$, even the invertible systems in this family therefore fall outside the scope of dimension theorems requiring strong irreducibility. The branch shears enter through the tilted projections $(t,z)\mapsto Pz+at$. Zooming into a cylinder changes the slope $a$. For each orthogonal projection $P$, Theorem~\ref{hom} gives convergence of normalized entropy uniformly over bounded sets of slopes. This uniformity lets us count graph cylinders at two scales as the slope changes along an orbit.

We can describe the exceptional parameters directly in terms of the interpolation data. For any fixed number of branches, every vector of vertical factors with nonzero coordinates and sum different from one gives equality for all choices of interpolation points and all linear images (Corollary~\ref{cor:generic}). Section~\ref{sec:exceptional} describes the exceptional affine subspace for fixed interpolation data. It gives dimension bounds at reduced rank, criteria for equality at rank one, and examples with distinct Hausdorff and Minkowski dimensions.

We use Barnsley's interpolation construction~\cite{ba}, with arbitrary branch widths. In the scalar full-branch case, each horizontal contraction maps the entire interval $[0,1]$ onto its branch interval. Clearing the denominators in the unequal shear ratios of~\cite[Equation~(2.4)]{brs} gives a nonzero generator in \eqref{rank}. In higher rank this condition requires the differences of shears obtained by composing pairs of active branches in opposite orders to span the residual space. The rank condition also makes sense when a denominator vanishes. B\'ar\'any, Rams and Simon also obtain dimension formulae for piecewise expanding affine repellers, with applications to interpolation functions~\cite{brs2020}. Section~\ref{masssec} gives counterexamples to parts (a) and (c) of the main theorem of Massopust~\cite[pp.~8--9]{mass}. For part (c), the printed statement and the proof predict different dimensions. Our example rules out both values.

For the Minkowski dimension of scalar affine interpolants see Hardin and Massopust~\cite{hm} and Ruan, Su and Yao~\cite{rsy}. Barnsley, Elton, Hardin and Massopust~\cite{behm} introduced hidden variable interpolants as coordinate projections. Hardin and Massopust~\cite[Theorem~6, p.~544]{hmproj} computed the Minkowski dimension of a coordinate graph in their model of hidden variable interpolation surfaces. For earlier Hausdorff-dimension criteria see Gibert and Massopust~\cite{gm}. The planar triangular setting is studied in~\cite{brs2}, and dimension bounds for vector-valued interpolants are given in~\cite{vpv}.

An early planar-valued example is Bara\'nski's complexification of the Weierstrass function~\cite{baranski2002}. In arbitrary target dimension, Ren and Shen~\cite{rs} study the series $W(x)=\sum_{k\ge0}\lambda^k\phi(b^kx)$ with real analytic $1$-periodic $\phi:\mathbb R\to\mathbb R^d$, integer $b\ge2$ and $1/b<\lambda<1$. If $q_{\rm Lip}$ is the largest dimension of a subspace onto which $W$ has Lipschitz projection, the Hausdorff and Minkowski dimensions of their graph both equal $\min\{\log_{1/\lambda} b,\,1+(d-q_{\rm Lip})(1+\log_b\lambda)\}$. For FIFs with equal widths and one common scalar satisfying $1/N<|s|<1$, setting $b=N$ and $\lambda=|s|$ gives the same algebraic formula with $r$ in the place of $d-q_{\rm Lip}$. Ren and Shen also use an inverse theorem for entropy and induction on the dimension of the target space. In our affine setting the horizontal locations of cylinders give separation for typical slopes, and the local entropy estimate in Lemma~\ref{frames} holds uniformly over the orthogonal complements that occur in the induction. For unequal widths and vertical factors, permutation blocks preserve the span of the branch shear differences (Lemma~\ref{perms}).

For bivariate interpolation surfaces, see Ma\l{}ysz~\cite{malysz} on Minkowski dimension and B\'ar\'any and Verma~\cite{bv} on the Hausdorff dimension of wedding cake surfaces over triangular domains. Here the domain is an interval and the interpolated values are vectors.

The proof of the upper bound in Section~\ref{upperbound} covers graph pieces at a stopping scale and includes branches with zero vertical factor. The lower bound uses three properties of the branch dynamics.
First, two distinct horizontal cylinders of length $k$ in a homogeneous subsystem have left endpoints at least $\rho^k$ apart. In a tilted projection this time difference is the coefficient of the slope. H\"older regularity of the first marginal of the stationary slope measure then gives separation of the projected translations outside a set of slopes of exponentially small measure (Lemma~\ref{sep}).

Second, scalar vertical multiplication preserves every linear subspace. Hochman's multidimensional inverse theorem~\cite{hoch} extends his one-dimensional entropy theorem~\cite{hoch1}. We handle its subspaces $V_i$ by induction on the target dimension $p$. Induction gives an entropy bound in $V_i^\perp$, and saturation gives one along $V_i$. In Lemma~\ref{frames} we obtain a lower bound with a single scale threshold for all complements. We condition on the dyadic cubes in the target space $\R^p$ before projecting.

Third, the permutations of a fixed multiset of symbols have common horizontal and vertical contractions. Interchanging the last two symbols gives a difference of shears independent of their common prefix. Lemma~\ref{perms} uses these differences to transfer the rank condition and the entropy of prescribed frequencies to homogeneous subsystems. The graph estimate in Section~\ref{permsec} uses the uniform entropy of tilted projections and counts surviving cylinders at the horizontal and vertical scales. Section~\ref{rangesection} then passes to an arbitrary linear image by transforming the interpolation data. The two geometric counterexamples in Section~\ref{masssec} can be read independently.

Theorems~\ref{main} and~\ref{hom} and Corollary~\ref{cor:generic} have been formalized in Lean~4. The formalization assumes Hochman's inverse theorem as stated in Theorem~\ref{litH}. This suffices for Theorem~\ref{main}, Corollary~\ref{cor:generic} and Theorem~\ref{hom}(1). The proof of Theorem~\ref{hom}(2) also assumes Feng's exact-dimensionality theorem. Part~(3) follows from part~(1) and a corollary of the Hochman--Shmerkin local entropy averages theorem. The Lean proofs contain no \texttt{sorry} placeholders. The Lean source files are available to referees on request.

\section{Statement and special cases}\label{statement}
Let $n\ge1$ and $N\ge2$. Fix partition points
$0=t_0<\cdots<t_N=1$ and interpolation values $y_i\in\R^n$. We call the pairs $(t_i,y_i)$ the interpolation points. For
$0\le i<N$ put $b_i=t_{i+1}-t_i$ and choose $|s_i|<1$.
The interval $I_i=[t_i,t_{i+1}]$ is the $i$th branch interval.
A \emph{branch} means the restriction of the self-referential
construction to $I_i$, with horizontal map
$T_i(t)=t_i+b_it$ and vertical scaling factor $s_i$. We call branch $i$ \emph{active} if $s_i\ne0$ and a \emph{zero branch} if $s_i=0$.
Define
\[
 \mathcal C_{y_0,y_N}=\{\varphi\in C([0,1],\R^n):\varphi(0)=y_0,\ \varphi(1)=y_N\}.
\]
The Read--Bajraktarevi\'c operator $\mathcal T$ on this complete space is
specified on every branch by
\begin{equation}\label{FIFeq}
 (\mathcal T\varphi)(t_i+b_it)
 =s_i \varphi(t)+(1-t)(y_i-s_i y_0)+t(y_{i+1}-s_i y_N).
\end{equation}
The branch expressions agree at the shared partition points and preserve
both endpoint values, so $\mathcal T$ maps $\mathcal C_{y_0,y_N}$
to itself. Moreover,
$\|\mathcal T\varphi-\mathcal T\psi\|_\infty\le
(\max_i|s_i|)\|\varphi-\psi\|_\infty$.
The Banach fixed-point theorem therefore gives a unique
continuous interpolant $f$, which satisfies \eqref{FIFeq} with
$\mathcal Tf=f$. This is the standard fractal interpolation construction
of Barnsley~\cite{ba}, also described in~\cite{massbook}.
For $N=1$ the interpolation equation yields the affine function
joining the two endpoints, so this degenerate case is excluded
from the notation below.

Write $G=\{(t,f(t)):0\le t\le1\}$ for the graph and let
$\ell_{\rm aff}(t)=(1-t)y_0+t y_N$ be the \emph{endpoint affine interpolant},
namely the straight line segment between the endpoint values.
Set $u_i=y_i-\ell_{\rm aff}(t_i)$,
$\Delta_i=u_{i+1}-u_i$, and $E=\operatorname{span}\{u_i\}$.
We call $h=f-\ell_{\rm aff}$ the \emph{normalized interpolant}. The corresponding bi-Lipschitz shear gives
\[
 h(t_i+b_it)=u_i+\Delta_i t+s_i h(t).
\]
Iteration from zero shows $h([0,1])\subset E$.
Among the normalized interpolation points $(t_i,u_i)$, the two endpoints differ by a vector that spans the time axis. Consequently,
$r=\dim E=\dim\operatorname{aff}\{(t_i,y_i)\}-1$.
The vectors $\Delta_i$ span $E$ and sum to zero.
If $r=0$, the graph is a straight segment.
We write $\underline{\dim}_{\B}$ and $\overline{\dim}_{\B}$ for the
lower and upper Minkowski (box-counting) dimensions.

\begin{definition}[The graph exponents]\label{def:D}
Let $J=\{i:s_i\ne0\}$ be the set of active branch indices, put $a_i=|s_i|$, and write $A(x)=\sum_{i\in J}a_i^x$. For every integer $j\ge0$, define $D_j$ as the unique $\beta\ge1$ such that
\begin{equation}\label{Ddef}
 \max\left\{\sum_{i=0}^{N-1} b_i^\beta,\ 
       \sum_{i\in J} a_i^{\min(\beta,j)} b_i^{\max(\beta-j,0)}\right\}=1 .
\end{equation}
Write $D=D_r$.
\end{definition}
We have $D_0=1$. For $j\ge1$, equation~\eqref{Ddef} gives
\[
 D_j=\begin{cases}
  1, & A(1)\le1,\\
  q, & A(1)>1\ge A(j),\\
  j+\delta, & A(j)>1.
 \end{cases}
\]
In the second case $q\in(1,j]$ is the root of $A(q)=1$.
In the third, $\delta\in(0,1)$ is the root of $\sum_{i\in J}a_i^j b_i^\delta=1$.
Indeed, the first sum in \eqref{Ddef} decreases strictly from one. For $j=0$ the second sum is at most the first. For $j\ge1$ it equals $A(\beta)$ up to $j$ and $\sum_{i\in J}a_i^jb_i^{\beta-j}$ thereafter. These expressions agree at $j$ and decrease strictly to zero when $J\ne\varnothing$. At $\beta=j+1$ the second sum is strictly less than $\sum_i b_i=1$. If $J$ is empty it vanishes identically.

The rank condition on the active branches is
\begin{equation}\label{rank}
 \operatorname{span}\{(b_k-s_k)\Delta_i-(b_i-s_i)\Delta_k:i,k\in J\}=E .
\end{equation}
\begin{theorem}[The graph and its linear images]\label{main}
Let $A$ and $D_j$ be as in Definition~\ref{def:D}.
\par\smallskip\noindent\textup{\textbf{(a) The graph.}}
Let $G$ be the graph of the interpolant and $r=\dim E$ its residual rank.
Then $\overline{\dim}_{\B}G\le D_r$. If $A(1)\le1$ or the active branches satisfy the rank condition~\eqref{rank}, then
\[
  \dim_{\rm H}G=\underline{\dim}_{\B}G=\overline{\dim}_{\B}G=D_r .
\]
\textup{\textbf{(b) Linear images.}}
Let $\Pi:\R\times\R^n\to\R^m$ be linear. Write $j$ for the dimension of the image of the residual space and $g$ for the affine dimension of the projected interpolation points, so that
\[
 j=\dim\Pi(\{0\}\times E),\qquad
 g=\dim\operatorname{aff}\{\Pi(t_i,y_i):0\le i\le N\}.
\]
The projected endpoint displacement either lies in the projected residual space or adds one independent direction. Thus $g\in\{j,j+1\}$ (Section~\ref{rangesection}). We have $\overline{\dim}_{\B}\Pi(G)\le\min\{D_j,g\}$. Under the same equality hypotheses as in \textup{(a)},
\begin{equation}\label{allrangeformula}
 \dim_{\rm H}\Pi(G)=\underline{\dim}_{\B}\Pi(G)=
 \overline{\dim}_{\B}\Pi(G)=\min\{D_j,g\}.
\end{equation}
If $s_i>0$ for every branch $i$, both equalities hold for every choice of interpolation values.
\end{theorem}
For a fixed map $\Pi(t,y)=tv+Ly$, let $W$ be the span on the left of \eqref{rank}. The equality in part~\textup{(b)} already holds if $L(W)=L(E)$: the transformed interpolation data have residual space $L(E)$ and their rank generators span $L(W)$ (Section~\ref{rangesection}). Thus an image may satisfy the formula even when $W$ is a proper subspace of $E$.
Taking $\Pi$ to be the identity in part~\textup{(b)} gives $j=r$ and $g=r+1$, hence part~\textup{(a)} since $D_r\le r+1$.

\begin{corollary}[Generic vertical factors]\label{cor:generic}
Fix $N\ge2$. For $s\in(-1,1)^N$ outside
\[
 \mathcal X_N=
 \bigcup_{i=0}^{N-1}\{s:s_i=0\}\ \cup\ \Bigl\{s:\sum_{i=0}^{N-1}s_i=1\Bigr\},
\]
the dimension equalities in Theorem~\ref{main} hold simultaneously for every choice of interpolation points with $N$ branches and every linear image of the graph.
This set of good parameters is open, dense and of full Lebesgue measure in $(-1,1)^N$.
\end{corollary}
\begin{proof}
All branches are active. Summing the generators in \eqref{rank} over $k$ gives
$(1-\sum_k s_k)\Delta_i$ because $\sum_k b_k=1$ and $\sum_k\Delta_k=0$.
The $\Delta_i$ span $E$ and the coefficient is nonzero, so Theorem~\ref{main} applies.
The excluded set is a finite union of proper affine hyperplanes intersected with the parameter cube.
\end{proof}
For positive factors, the remaining case $\sum_i s_i=1$ gives $A(1)=1$. Together with Corollary~\ref{cor:generic}, this proves the final assertion of Theorem~\ref{main}.

Under the equality hypotheses, taking $\Pi(t,y)=(t,Ly)$ gives dimension $D_{\dim L(E)}$ for the graph of $Lf$. This includes the graph of each scalar coordinate. Further range formulae follow in Section~\ref{rangesection}.

Only the moduli $a_i$ enter $D$. The signs of the $s_i$ enter through \eqref{rank}.
For equal widths $1/N$ and a common scalar $s$, assume $r>0$ and $1/N<|s|<1$. Setting $\alpha=\log(1/|s|)/\log N$ and $q=1/\alpha$ gives
\[
 D=\min\{1/\alpha,\,1+r(1-\alpha)\}.
\]
Indeed, the second candidate minus the first is $(1-\alpha)(r-q)$, so the transition is exactly at $r=q$. Equality of dimensions follows from Theorem~\ref{main} also for a common negative scalar, since then $\sum_i s_i\ne1$. For $|s|\le1/N$, or for $r=0$, the dimension is one. In particular the canonical data $0,e_1,\ldots,e_{N-1},0$ have rank $N-1$. Table~\ref{tab:threebranch} gives the dimension for three equal branches and four noncoplanar interpolation points in $\R^3$.

\begin{table}[htbp]
\centering
\caption{Graph dimension for three equal branches, residual rank two and a common vertical factor $s>0$.}\label{tab:threebranch}
\smallskip
\begin{tabular}{lc}\toprule
Parameter & $\dim_{\mathrm H}G=\dim_{\mathrm B}G$\\\midrule
$0<s\le1/3$ & $1$\\
$1/3<s<1/\sqrt3$ & $\log3/\log(1/s)$\\
$s=1/\sqrt3$ & $2$\\
$1/\sqrt3<s<1$ & $3+2\log_3s$\\\bottomrule
\end{tabular}
\end{table}

At $s=1/\sqrt2$ this gives $3-\log_3 2$. For unequal positive vertical factors $(1/2,2/3,3/4)$ and the same partition and rank, the answer is $2+\log_3(181/144)$.

\section{Upper bound}\label{upperbound}
\begin{definition}[Normalized affine system]\label{def:corrsys}
Identify $E$ with $\R^r$ and put $d=r+1$.
In this section $\beta$ is an auxiliary dimension exponent.
The symbol $t$ remains reserved for the horizontal parameter of the FIF.
In the coordinate order $(z,t)$, define the normalized affine branch maps by
\[
 F_i(z,t)=(s_i z+u_i+\Delta_i t,\ t_i+b_i t),
 \qquad i=0,\ldots,N-1.
\]
Their linear parts are
\[
 L_i=\begin{pmatrix}s_i I_r&\Delta_i\\0&b_i\end{pmatrix}.
\]
\end{definition}
The upper-right block of the commutator $L_iL_j-L_jL_i$
is $(b_j-s_j)\Delta_i-(b_i-s_i)\Delta_j$, the vector appearing
in the rank condition~\eqref{rank}.
Write $\widetilde G=\{(h(t),t):0\le t\le1\}$ for the normalized graph.
Then $\widetilde G=\bigcup_i F_i(\widetilde G)$.
The affine coordinate change $(t,f(t))\mapsto(h(t),t)$ is
bi-Lipschitz, so $G$ and $\widetilde G$ have the same Hausdorff and
Minkowski dimensions.

In this section a \emph{word} is a finite sequence $w=(i_1,\ldots,i_k)$ with entries in
$\{0,\ldots,N-1\}$, of any length $k\ge0$. Words may contain zero branches. For a word $w$ put
\[
 T_{w}=T_{i_1}\circ\cdots\circ T_{i_k},\qquad
 b_{w}=b_{i_1}\cdots b_{i_k},\qquad
 s_{w}=s_{i_1}\cdots s_{i_k},\qquad
 a_{w}=|s_{w}| ,
\]
so that $T_{w}(t)=T_{w}(0)+b_{w}t$. For the empty word, $T_\varnothing$ is the identity and
$b_\varnothing=s_\varnothing=1$. The word $w$ followed by the index $i$ is written $w i$.
We continue to identify $E$ with $\R^r$. We use the Euclidean norm $\|\cdot\|$ and take all cubes to be axis-parallel in
$\R^r\times\R=\R^{d}$. Write $N_\varepsilon(\mathcal K)$ for the minimum number of axis-parallel cubes of side $\varepsilon$ covering a bounded set $\mathcal K$.

\begin{lemma}[Pieces of the normalized graph]\label{pieces}
Define vectors $u_{w},\Delta_{w}\in E$ for all words by $u_\varnothing=\Delta_\varnothing=0$ and
\begin{equation}\label{piecerec}
 u_{w i}=u_{w}+t_i\Delta_{w}+s_{w}u_i,\qquad
 \Delta_{w i}=b_i\Delta_{w}+s_{w}\Delta_i .
\end{equation}
For a word of length one these are the vectors $u_i$ and $\Delta_i$ of Section~\ref{statement}. For every word $w$,
\begin{equation}\label{pieceeq}
 h(T_{w}(t))=u_{w}+\Delta_{w}t+s_{w}h(t)\qquad(0\le t\le1).
\end{equation}
Put $\widetilde G_{w}=\{(h(T_{w}(t)),T_{w}(t)):0\le t\le1\}$. Then $\widetilde G_\varnothing=\widetilde G$,
\begin{equation}\label{piecesplit}
 \widetilde G_{w}=\bigcup_{i=0}^{N-1}\widetilde G_{w i},
\end{equation}
and
\begin{equation}\label{piecechild}
 a_{w i}=a_{w}a_i,\qquad b_{w i}=b_{w}b_i,\qquad
 \|\Delta_{w i}\|\le b_i\|\Delta_{w}\|+a_{w}\|\Delta_i\| .
\end{equation}
\end{lemma}
\begin{proof}
Equation~\eqref{pieceeq} holds for the empty word. Substituting $t_i+b_it$ for $t$ in the identity for $w$ and applying the branch equation for $h$ gives the coefficients in \eqref{piecerec}. This proves \eqref{pieceeq} by induction. The branch intervals cover $[0,1]$, giving \eqref{piecesplit}. The product definitions and the triangle inequality in \eqref{piecerec} give \eqref{piecechild}.
\end{proof}

\begin{lemma}[Two coverings of one piece]\label{twocovers}
Let $C_h=\max_{0\le t\le1}\|h(t)\|$, let $w$ be a word and let $\varepsilon>0$.
\begin{enumerate}
\item \emph{Thin tube.} If $a_{w}C_h\le\varepsilon$, then
$N_{4\varepsilon}(\widetilde G_{w})\le(b_{w}+\|\Delta_{w}\|)/\varepsilon+2$.
\item \emph{Narrow column.} If $b_{w}\le\varepsilon$, then
$N_{\varepsilon}(\widetilde G_{w})\le\big(2(a_{w}C_h+\|\Delta_{w}\|)/\varepsilon+1\big)^r$.
\end{enumerate}
\end{lemma}
\begin{proof}
(1) By \eqref{pieceeq} the point of $\widetilde G_{w}$ with parameter $t$ is within Euclidean
distance $a_{w}\|h(t)\|\le\varepsilon$ of the point $(u_{w}+\Delta_{w}t,\ T_{w}(t))$. These points
form a segment of length at most $b_{w}+\|\Delta_{w}\|$. Let $n_{\rm seg}\ge1$ be the least integer with
$n_{\rm seg}\varepsilon\ge b_{w}+\|\Delta_{w}\|$. Every point of the segment is within distance
$\varepsilon$ of one of the $n_{\rm seg}+1$ points with parameters $j/n_{\rm seg}$, $0\le j\le n_{\rm seg}$. Hence
every point of $\widetilde G_{w}$ is within distance $2\varepsilon$ of one of these $n_{\rm seg}+1$ points
and lies in the cube of side $4\varepsilon$ centred there. Finally
$n_{\rm seg}+1\le(b_{w}+\|\Delta_{w}\|)/\varepsilon+2$.

(2) By \eqref{pieceeq} and the inclusion $h([0,1])\subset E$, the first coordinate of every point of
$\widetilde G_{w}$ differs from $u_{w}$ by a vector of $E=\R^r$ of norm at most
$a_{w}C_h+\|\Delta_{w}\|$. It therefore lies in a cube of $\R^r$ of side
$2(a_{w}C_h+\|\Delta_{w}\|)$, and such a cube is covered by at most
$\big(2(a_{w}C_h+\|\Delta_{w}\|)/\varepsilon+1\big)^r$ cubes of side $\varepsilon$. The second
coordinate lies in $T_{w}([0,1])$, an interval of length $b_{w}\le\varepsilon$. Hence
$\widetilde G_{w}$ is covered by that many cubes of side $\varepsilon$ in $\R^r\times\R$.
\end{proof}

\begin{proposition}[Direct covering bound]\label{direct}
Let $r\ge1$ and let $\beta$ be a real number satisfying one of the following conditions.
\begin{enumerate}
\item[(A)] $\beta\ge1$ and $A(\beta)<1$.
\item[(B)] $\beta>r$ and $B_r(\beta)<1$, where $B_r(x)=\sum_{i\in J}a_i^rb_i^{x-r}$.
\end{enumerate}
Then there is a constant $C_{\rm dir}$, depending only on the interpolation data and on $\beta$,
such that
\[
 N_{4\varepsilon}(\widetilde G)\le C_{\rm dir}\,\varepsilon^{-\beta}\qquad(0<\varepsilon\le1).
\]
In particular $\overline{\dim}_{\B}G\le\beta$. The set $J$ is arbitrary: zero branches and the
case $J=\varnothing$ are included.
\end{proposition}
\begin{remark}\label{directrem}
In case (A) a word
is stopped when the vertical size $a_{w}C_h$ of its piece falls below $\varepsilon$, and the
piece is then a thin tube around a segment. In case (B) a word is stopped when its width
$b_{w}$ falls below $\varepsilon$, and the piece is then a narrow column over a ball of $E$. A
zero branch gives $a_{w i}=0$, so its piece is thin for every $\varepsilon$ and needs no
separate treatment.
\end{remark}

\begin{proof}
Write $C_h=\max\|h(t)\|>0$ and $T_\Delta=\sum_i\|\Delta_i\|$. The quantities in the two stopping rules may be summed over all branches, since $a_i=0$ off $J$:
\begin{equation}\label{directsums}
 \sum_i a_i^\beta=A(\beta),\qquad
 \sum_i a_i^r b_i^{\beta-r}=B_r(\beta).
\end{equation}
In each case take $0<\varepsilon\le1$. By \eqref{piecesplit}, the graph pieces corresponding to the leaves cover the normalized graph. Their horizontal widths $b_w$ sum to one.

\emph{Case (A).} Put $\tau=\varepsilon/C_h$ and stop each word the first time $a_w\le\tau$. The tree is finite since $\max_i a_i<1$. Write $\mathcal I$ for its internal vertices and $\mathcal L$ for its leaves. Summing $a_w^\beta$ over all words gives $\sum_{k\ge0}A(\beta)^k=(1-A(\beta))^{-1}$. Since $a_w>\tau$ on $\mathcal I$ and $\beta\ge1$,
\[
 |\mathcal I|\le\frac{\tau^{-\beta}}{1-A(\beta)},\qquad
 \sum_{w\in\mathcal I}a_w\le\frac{\tau^{1-\beta}}{1-A(\beta)}.
\]
Every internal vertex has $N$ children, so $|\mathcal L|=1+(N-1)|\mathcal I|$. By \eqref{piecechild},
\begin{equation}\label{directlen}
 \sum_i\|\Delta_{wi}\|\le\|\Delta_w\|+a_w T_\Delta .
\end{equation}
Summing \eqref{directlen} over the internal vertices cancels all internal shear terms, since $\Delta_\varnothing=0$. Thus $\sum_{w\in\mathcal L}\|\Delta_w\|\le T_\Delta\sum_{w\in\mathcal I}a_w$. On every leaf Lemma~\ref{twocovers}(1) applies, and hence
\[
 N_{4\varepsilon}(\widetilde G)
 \le 2|\mathcal L|+\varepsilon^{-1}\Big(1+T_\Delta\sum_{w\in\mathcal I}a_w\Big)
 =O(\varepsilon^{-\beta}).
\]
This also covers the tree consisting only of its root.

\emph{Case (B).} Put $q_i=b_i^{\beta-r}$ and define
\[
 \begin{aligned}
 A_b&=\Big(\sum_i b_i^\beta\Big)^{1/r}<1,&
 B_a&=B_r(\beta)^{1/r}<1,\\
 C_\Delta&=\Big(\sum_i q_i\|\Delta_i\|^r\Big)^{1/r},&
 K&=\max\{C_h,C_\Delta/(1-B_a)\}.
 \end{aligned}
\]
Here $A_b<1$ since $\beta>r\ge1$ and $\sum_i b_i=1$. For $d_w=\|\Delta_w\|+K a_w$, apply Minkowski's inequality with weights $q_i$ to \eqref{piecechild}. It gives
\[
 \Big(\sum_i q_i d_{wi}^{\,r}\Big)^{1/r}
 \le A_b\|\Delta_w\|+(C_\Delta+K B_a)a_w\le d_w .
\]
Consequently $V(w)=b_w^{\beta-r}d_w^r$ satisfies $\sum_i V(wi)\le V(w)$. Stop at the first $b_w\le\varepsilon$. This tree is finite since $\max_i b_i<1$. On its leaves $b_{\min}\varepsilon\le b_w\le\varepsilon$, where $b_{\min}=\min_i b_i>0$. Their widths sum to one and their $V$-values sum to at most $V(\varnothing)=K^r$, so
\[
 |\mathcal L|\le (b_{\min}\varepsilon)^{-1},\qquad
 \sum_{w\in\mathcal L}d_w^r\le (b_{\min}\varepsilon)^{-(\beta-r)}K^r .
\]
Lemma~\ref{twocovers}(2), $K\ge C_h$ and $(1+2x)^r\le 2^{r-1}(1+2^r x^r)$ for $x\ge0$ now give
\[
 N_{4\varepsilon}(\widetilde G)\le N_\varepsilon(\widetilde G)
 \le 2^{r-1}\left(|\mathcal L|+2^r\varepsilon^{-r}\sum_{w\in\mathcal L}d_w^r\right)
 =O(\varepsilon^{-\beta}).
\]
All implied constants depend only on the interpolation data and $\beta$. The bi-Lipschitz coordinate change transfers the upper Minkowski dimension bound to $G$.
\end{proof}

\begin{proof}[Graph upper bound]
Let $\beta>D$. Suppose first that $A(r)\le1$. Since $A(r)\le A(1)$, this covers the first two
cases of \eqref{Ddef}. In the first case $D=1$ and $A(D)=A(1)\le1$. In the second case
$A(D)=A(q)=1$. If $J=\varnothing$, then $A(\beta)=0$. Otherwise $A$ is strictly decreasing,
so $A(\beta)<A(D)\le1$. Since also $\beta>D\ge1$, condition (A) of Proposition~\ref{direct}
holds. Suppose now that
$A(r)>1$. Then $J\ne\varnothing$, $D=r+\delta$ and
$B_r(D)=\sum_{i\in J}a_i^rb_i^\delta=1$. As $B_r$ is strictly decreasing, $B_r(\beta)<1$,
and $\beta>r$, so condition (B) holds. In either case Proposition~\ref{direct} gives
$\overline{\dim}_{\B}G\le\beta$. Letting $\beta\downarrow D$ yields
$\overline{\dim}_{\B}G\le D$.
\end{proof}

For the singular-value pressure of the active triangular matrices, see
Falconer and Miao~\cite{fm}.

\section{A homogeneous projection theorem}\label{sec:homogeneous}
For the lower bound we use an auxiliary homogeneous system. Its maps
$\mathcal S_i$ are distinct from the normalized FIF maps $F_i$
of Section~\ref{upperbound}.
Consider $M\ge2$ maps
\[
 \mathcal S_i(t,z)=(x_i+\rho t,d_i+v_it+Sz),\qquad z\in\R^r,\quad 0<\rho<S<1.
\]
We write the time coordinate first, reversing the order used in Section~\ref{upperbound}. The linear parts are
\[
 \mathcal L_i=\begin{pmatrix}\rho&0\\v_i&SI_r\end{pmatrix},
 \qquad \mathcal L_w=\mathcal L_{i_1}\cdots\mathcal L_{i_k}
 \quad\text{for }w=(i_1,\ldots,i_k).
\]
Fix a probability vector $\boldsymbol\pi=(\pi_1,\ldots,\pi_M)$ with $\pi_i>0$ and $\sum_i\pi_i=1$. Let $\mu_{\boldsymbol\pi}$ be the self-affine measure obtained as the push-forward of the Bernoulli measure with weights $\boldsymbol\pi$ under the coding map, and let $(T,Z)\sim\mu_{\boldsymbol\pi}$. In this section logarithms have base two. Put
\[
 \ell_S=\log(1/S),\quad X=\log(1/\rho),\quad
 \hw=H(\boldsymbol\pi)=-\sum_i\pi_i\log\pi_i,\quad Q=\hw/\ell_S.
\]
For $1\le p\le r$ and $P:\R^r\to\R^p$ with orthonormal rows let $\nu_a$ be the push-forward of $\mu_{\boldsymbol\pi}$ under $(t,z)\mapsto Pz+at$, where $a\in\R^p$. Write $H_\ell$ for unnormalized entropy of half-open dyadic cells of side $2^{-\ell}$.
\begin{theorem}\label{hom}
Let $\mathcal S_i(t,z)=(x_i+\rho t,d_i+v_i t+Sz)$, $1\le i\le M$, with $M\ge2$ and $0<\rho<S<1$. Suppose the intervals $x_i+\rho[0,1]\subset[0,1]$ have disjoint interiors and the differences $v_i-v_{i'}$ span $\R^r$. Let $\boldsymbol\pi$ be any strictly positive probability vector on these branches.
\begin{enumerate}
\item For every $1\le p\le r$, every $P:\R^r\to\R^p$ with orthonormal rows and every $B_0>0$,
\[
 \lim_{\ell\to\infty}\ \sup_{\|a\|\le B_0}\left|\frac{H_\ell(\nu_a)}{\ell}-\min(p,Q)\right|=0 .
\]
\item The measure $\mu_{\boldsymbol\pi}$ is exact dimensional, of dimension
\[
 \hw/X+(1-\ell_S/X)\min(r,Q).
\]
\item For each $1\le p\le r$, each $P$ as in part (1) and each fixed $a\in\R^p$,
the measure $\nu_a$ is exact dimensional with dimension $\min(p,Q)$.
\end{enumerate}
\end{theorem}
The projected branch images may overlap. The value in part (2) is the Lyapunov dimension for weights $\boldsymbol\pi$, with positive contraction exponents $\ell_S$ (of multiplicity $r$) and $X$.
Parts (2) and (3) are proved in Appendix~\ref{measureappendix}. If some weights vanish, discard their branches and verify the shear span for the remaining alphabet.
The measure formulae concern homogeneous systems with common $\rho,S$. For the Ledrappier--Young
formula and related exact-dimensionality results see \cite{ly,bk,feng}.
\paragraph{Plan of the proof of part (1).} Fix $P$. Lemma~\ref{loc} constructs the entropy rate $\kappa(P)$ averaged over the stationary distribution of slopes and gives a uniform lower bound in terms of this rate. If $\kappa(P)<\min(p,Q)$, the separation estimate of Lemma~\ref{sep} lets us extract a factor of positive entropy whose convolution with a tilted projection has negligible entropy growth (Lemma~\ref{extraction}). Hochman's inverse theorem then supplies subspaces at the component scales. Proposition~\ref{kappaprop} uses the lower-rank projections onto their orthogonal complements to rule out the deficit. The uniform upper bound completes the proof.
\subsection{Entropy estimates used in the proof}
Throughout this section $H$ is Shannon entropy in bits, without division by scale. For a probability measure $\omega$ on $\R^d$, $\mathcal D_i(y)$ denotes the half-open dyadic cube containing $y$ and $H_i(\omega)=H(\mathcal D_i(Y))$ for $Y\sim\omega$. Write $h_m(\omega)=H_m(\omega)/m$ for $m\ge1$. If $C=\mathcal D_i(y)$ has positive mass and lower corner $c$, then $\omega^{y,i}$ is the conditional measure on $C$ pushed forward by $u\mapsto2^i(u-c)$. When we average over $y$, each cell $C$ receives weight $\omega(C)$.

\begin{lemma}[Grid comparisons]\label{grid}
Let $I,J'$ be random variables with finitely many values on one probability space. If every value of $I$ occurs with positive probability together with at most $N_{\rm cmp}$ values of $J'$, and every value of $J'$ together with at most $N_{\rm cmp}$ values of $I$, then $|H(I)-H(J')|\le\log N_{\rm cmp}$. In particular, with
\[
 C_{\rm grid}(d,B)=d\log(\lceil B\sqrt d\rceil+2)\qquad(d\ge1,\ B\ge1),
\]
the following hold.
\begin{enumerate}
\item If $Y,Y'$ are bounded random vectors in $\R^d$ on one probability space and $\|Y-Y'\|_\infty\le A2^{-m}$ almost surely, the level-$m$ entropies of their distributions differ by at most $d\log(2\lceil A\rceil+3)$.
\item Let $\omega$ be a probability measure with bounded support in $\R^{d_1}$, let $\Psi$ map its support into $\R^{d_2}$, where $d_1,d_2\le d$, and let $m,m'\ge0$ be integers with
\[
 B^{-1}2^{m}\|y-y'\|\le2^{m'}\|\Psi(y)-\Psi(y')\|\le B2^{m}\|y-y'\|
\]
for all $y,y'$ in the support. Then $|H_{m'}(\Psi_*\omega)-H_m(\omega)|\le C_{\rm grid}(d,B)$. This covers a change to a translated grid, an isometry onto a Euclidean subspace with the grid of orthonormal coordinates on that subspace, and a change of the mesh by a factor in $[1/B,B]$.
\item For probability measures $\xi,\omega$ with bounded supports in $\R^d$, $H_m(\xi*\omega)\ge H_m(\omega)-C_{\rm grid}(d,1)$.
\end{enumerate}
\end{lemma}
\begin{proof}
For finite-valued random variables $I,J'$, use $H(I)-H(J')=H(I\mid J')-H(J'\mid I)$. The conditional entropies are bounded by the logarithms of the respective numbers of possible partners. In (1) the coordinate indices differ by at most $\lceil A\rceil+1$. For (2) the part of the support in any source cube of side $2^{-m}$ has diameter at most $\sqrt d\,2^{-m}$ and maps to a set of diameter at most $B\sqrt d$ in units of the target mesh $2^{-m'}$. It meets at most $(\lceil B\sqrt d\rceil+2)^d$ target cells. Applying the same argument to $\Psi^{-1}$ on its image gives the reverse inequality. This also covers grids in different orthonormal coordinates on one subspace. For (3), condition on the summand with distribution $\xi$ and use (2) for a translate.
\end{proof}

\begin{lemma}[Entropy increments]\label{increments}
For integers $m,K\ge1$ and compactly supported $\omega$ in $\R^d$,
\begin{align}
 \E_y H_m(\omega^{y,i})&=H_{i+m}(\omega)-H_i(\omega),\label{inc}\\
 \left|\frac1K\sum_{i=0}^{K-1}\E_y h_m(\omega^{y,i})
 -\frac{H_K(\omega)-H_0(\omega)}K\right|&\le dm/K.\label{tele}
\end{align}
If $\omega$ is a mixture with parameter $W$, and each conditional support meets at most $N_{\rm cell}$ level-$i$ cells, then
\begin{equation}\label{mix}
 H_{i+m}(\omega)-H_i(\omega)\ge H(\mathcal D_{i+m}(Y)\mid W)-\log N_{\rm cell}.
\end{equation}
\end{lemma}
\begin{proof}
Dyadic nesting identifies the left side of \eqref{inc} with entropy of the fine index given the coarse one. In \eqref{tele}, expand each increment into $m$ one-step increments. The sum is the average of the $m$ windows $H_{K+t}-H_t$, $0\le t<m$, divided by $K$. Each window differs from $H_K-H_0$ by at most $dt$ since a one-step increment belongs to $[0,d]$. For \eqref{mix}, further conditioning the fine index on $W$ can only reduce its conditional entropy. The chain rule for the nested indices then subtracts at most $H(\mathcal D_i(Y)\mid W)\le\log N_{\rm cell}$.
\end{proof}

\subsection{Cylinders, slopes and uniform local bounds}
For a prefix $w$ in outer-to-inner order write $x_w,d_w$ for its horizontal and vertical translations. Put $\lambda=\rho/S$ and $R_i a=\lambda a+Pv_i/S$. Applying $Pz+at$ to a single branch yields the constant $Pd_i+ax_i$ and the tail $S(Pz+(R_i a)t)$. Iterating gives the conditional projection
\begin{equation}\label{cyl}
 c_w(a)+S^{|w|}Y_w,\qquad c_w(a)=Pd_w+ax_w,\qquad
 Y_w\sim\nu_{R_w a},\quad R_w=R_{w_{|w|}}\circ\cdots\circ R_{w_1}.
\end{equation}
Fix a slope radius $B_0$, at least the radius in Theorem~\ref{hom}(1), satisfying $B_0\ge\max_i\|v_i\|/[S(1-\lambda)]$. In each target dimension the ball $B$ of this radius is invariant. The stationary slope measure $\eta=\eta_{P,\boldsymbol\pi}$ is the self-similar measure for the contractions $R_i$ with weights $\pi_i$. It is the distribution of the convergent random series of their translations with coefficients $\lambda^j$, $j\ge0$, and independent symbols of weights $\boldsymbol\pi$. If $a_*\sim\eta$ is independent of a $\boldsymbol\pi$-distributed word $w$, then $R_w a_*\sim\eta$ and $\|R_w a-R_w a_*\|\le2B_0\lambda^{|w|}$ for $a\in B$.

\paragraph{Dependence of constants.}
Fix the homogeneous system and $B_0$. Since $0\le T\le1$ and $Z$ is bounded, all $PZ+aT$ lie in a ball of radius $B_Y=\sup\|Z\|+B_0$, uniformly for every orthonormal-row map $P$. The grid comparison bounds below follow from Lemma~\ref{grid} with dimension at most $2r$, fixed scale factors, and these support bounds. The constant $C_{\rm cyl}$ bounds the error from rescaling a cylinder and conditioning on its coarse index. The constant $C_{\rm loc}$ also includes the errors from coupling slopes and comparing projected grids. Both bounds are uniform in $P,U,a$ and in all scale indices. The constants in Lemmas~\ref{sep} and~\ref{extraction} may depend on $P$, which is fixed in both applications. Thresholds declared after fixing $m$ or $n$ may depend on that parameter.

\begin{lemma}[Local entropy and all target directions]\label{loc}\label{frames}
For an orthonormal-row $P:\R^r\to\R^p$, put $k_i=\lceil i/\ell_S\rceil$ and
$\mathcal E_m(P)=\int H_m(\nu_a)\,d\eta_P(a)$. The limit
$\kappa(P)=\lim_m\mathcal E_m(P)/m$ exists and is at most $\min(p,Q)$.
Put $\kappa_0(P)=0$ and $\kappa_e(P)=\inf_{UU^*=I_e}\kappa(UP)$ for
$U:\R^p\to\R^e$, $1\le e\le p$. For every $\gamma>0$ there is $m_0$ such
that for every $m\ge m_0$ there is $i_0(m)$ with
\begin{equation}\label{comps}
 \E_y h_m(U_*\nu_a^{y,i})\ge\kappa_e(P)-\gamma
 \qquad(0\le e\le p,\ UU^*=I_e,\ a\in B,\ i\ge i_0(m)).
\end{equation}
Here conditioning is in $\R^p$ before projection. Moreover
$\liminf_n\inf_{a\in B} h_n(\nu_a)\ge\kappa(P)$.
\end{lemma}
\begin{proof}
\emph{Stationary entropy rate.} The ratio $S^{k_i}/2^{-i}$ lies in $(S,1]$. Every expectation over words in this section uses the product weights $\pi_w=\prod_{j}\pi_{w_j}$. Each conditional support in \eqref{cyl} has diameter at most $2B_Y2^{-i}$, hence meets at most $(\lceil2B_Y\rceil+2)^r$ coarse cells. Given the prefix, its fine entropy differs from $H_m(\nu_{R_wa})$ by at most $C_{\rm grid}(r,1/S)$. Applying \eqref{mix} gives
\begin{equation}\label{super}
 H_{i+m}(\nu_a)-H_i(\nu_a)\ge \E_wH_m(\nu_{R_w a})-C_{\rm cyl}.
\end{equation}
The constant $C_{\rm cyl}$ depends only on the fixed system and slope ball. Integrating \eqref{super} against the stationary slope measure shows $\mathcal E_n(P)-C_{\rm cyl}$ is superadditive. A bounded support meets $O(2^{pi})$ cells. Writing $W$ for the weighted prefix of length $k_i$, we have $H(W)=k_i\hw$. Moreover $H_i(\nu_a)\le H(W)+O(1)$, since the conditional support within each cylinder meets only boundedly many scale-$i$ cells. Thus $H_i(\nu_a)\le i\min(p,Q)+O(1)$ uniformly. Fekete's lemma applied to $\mathcal E_m(P)-C_{\rm cyl}$ now gives the limit.

\emph{Uniformity over projected directions.} The rank-zero projection is constant and has zero entropy. Fix $1\le e\le p$. Temporarily write $R_i^P,\nu_a^P$ to indicate the projection. Direct substitution gives $UR_i^P=R_i^{UP}U$ and $U_*\nu_a^P=\nu_{Ua}^{UP}$. The stationary series therefore gives $U_*\eta_P=\eta_{UP}$, and
\[
 \mathcal E_\ell(UP)=\int H_\ell(U_*\nu_a^P)\,d\eta_P(a).
\]
At a fixed orthonormal frame $U_0$, the definition of $\kappa_e(P)$ and the limit just proved allow a length $\ell\ge1$ with
\[
 \mathcal E_\ell(U_0P)\ge \ell(\kappa_e(P)-\gamma/4)+C_{\rm cyl}+r\log5.
\]
If $\|U-U_0\|_{\rm op}<2^{-\ell}/(1+B_Y)$, evaluating $U$ and $U_0$ on the same point of $\nu_a^P$ gives a displacement at most $2^{-\ell}$. Lemma~\ref{grid}(1) implies $\mathcal E_\ell(UP)-C_{\rm cyl}\ge \ell(\kappa_e(P)-\gamma/4)$ throughout this neighborhood. The set $UU^*=I_e$ is closed and bounded in a finite-dimensional matrix space. Take a finite subcover and let $\ell_*$ be the largest of its lengths. For any $U$, select a neighborhood containing it and denote the associated length by $\ell$. Write $m=j \ell+s$, $0\le s<\ell$. Iterating the stationary version of \eqref{super} yields
\[
 \mathcal E_m(UP)\ge j(\mathcal E_\ell(UP)-C_{\rm cyl})
 \ge m(\kappa_e(P)-\gamma/4)-p \ell_* .
\]
For $s>0$ the discarded term is $\mathcal E_s(UP)\ge0$. For $s=0$ at most $j-1$ subtractions occur. Hence $\mathcal E_m(UP)/m\ge \kappa_e(P)-\gamma/2$ when $m\ge4p \ell_*/\gamma$. Take the maximum over the finitely many positive ranks $e\le p$.

\emph{Conditioning at an arbitrary slope.} It remains to connect this estimate to the conditioning in \eqref{comps}. Let $Y\sim\nu_a^P$, $W$ its random prefix of length $k_i$, $C_i=\mathcal D_i(Y)$ and $J_{i+m}=\mathcal D_{i+m}(UY)$. The conditional entropy chain rule gives
\[
 H(J_{i+m}\mid C_i)\ge H(J_{i+m}\mid W,C_i)
 \ge H(J_{i+m}\mid W)-H(C_i\mid W).
\]
The last subtracted term is at most $r\log(\lceil2B_Y\rceil+2)$. By \eqref{cyl} and Lemma~\ref{grid}(2), the first term on the right differs from $\E_wH_m(U_*\nu_{R_wa}^P)$ by at most $C_{\rm grid}(r,1/S)$. Given $C_i$, normalizing its lower corner translates the projected grid, so $\E_y H_m(U_*\nu_a^{y,i})$ differs from $H(J_{i+m}\mid C_i)$ by at most $C_{\rm grid}(r,1)$. Choose $i_0(m)$ so that $2B_0\lambda^{k_i}\le2^{-m}$ for $i\ge i_0(m)$. Iterate the slope maps from $a$ and from a stationary initial slope, using the same word and the same terminal coding point in the two projections. Their terminal slopes differ by at most $2B_0\lambda^{k_i}$. Since $0\le T\le1$ and $\|U\|=1$, the resulting grid entropies differ by at most $r\log5$. Combining the bounds gives
\[
 \E_y H_m(U_*\nu_a^{y,i})\ge \mathcal E_m(UP)-C_{\rm loc},
 \qquad C_{\rm loc}=C_{\rm cyl}+C_{\rm grid}(r,1)+r\log5 .
\]
Taking also $m\ge2C_{\rm loc}/\gamma$ proves the assertion. The thresholds $m$ and $i_0(m)$ are uniform over all choices of the frames $U_i$.
For $e=p$ an orthogonal $U$ changes grid entropy by $O_r(1)$ at every level,
so $\kappa(UP)=\kappa(P)$. Taking $U=I_p$ and summing consecutive $m$-block
increments from $i_0(m)$ proves the final assertion after $n\to\infty$ and
$\gamma\downarrow0$. The case $p=0$ has zero entropy.
\end{proof}

\subsection{Separation from the time coordinate}
\begin{lemma}[H\"older regularity of the first slope marginal and separation of translations]\label{sep}
Fix $P$ of rank $p\ge1$. There are constants $C_P<\infty$ and $\theta_P>0$ such that the first marginal of $\eta$ assigns mass at most $C_P t^{\theta_P}$ to every interval of length $t\le1$. There are an integer $d_{\rm sep}\ge1$ and constants $b_P>0$, $C_P'<\infty$ such that for every $k\ge1$ the following holds for all slopes $a$ outside a set of $\eta$-mass at most $C_P'2^{-b_Pk}$: the translations $c_w(a)$ for distinct words of length $k$ are pairwise more than $2^{-d_{\rm sep} k}$ apart in maximum norm.
\end{lemma}
\begin{proof}
The first row of $P$ is nonzero. Since the differences $v_i-v_j$ span, their projections onto that row cannot all be zero. The scalar slope maps thus have at least two distinct translations. The convex hull of their stationary support is an interval of length $D_*>0$. The extreme translations fix its endpoints, as follows by taking extrema in the support invariance equation. Choose $b$ with $\lambda^b<1/2$. The two $b$-fold extreme repetitions map the interval to opposite ends with gap $g=D_*(1-2\lambda^b)>0$. Iterate the stationarity equation in blocks of $b$ symbols. Within each parent of depth $bj$ the two extreme descendants have a gap $g\lambda^{bj}$, and each has relative weight at least $\pi_{\min}^b$, where $\pi_{\min}=\min_i\pi_i>0$. An interval of length $t<g\lambda^{b(n-1)}$ misses at least one child of conditional weight at least $\pi_{\min}^b$ at each surviving prefix in the first $n$ block levels. Induction on levels bounds its mass by $(1-\pi_{\min}^b)^n$. Choosing the largest such $n$ gives the power bound with $\theta_P=\log(1-\pi_{\min}^b)/(b\log\lambda)>0$. A fixed prefactor covers larger $t$ and endpoints of scale ranges.

At level $k$ the horizontal intervals have equal length $\rho^k$ and disjoint interiors: after the first differing symbol this follows by applying the common prefix to the first-level intervals. Thus $|x_w-x_{w'}|\ge\rho^k$. In the first coordinate, $c_w(a)-c_{w'}(a)$ has the form $u+(x_w-x_{w'})a_1$. The probability that its absolute value is at most $2^{-dk}$ is bounded by $C_P(2^{1-dk}\rho^{-k})^{\theta_P}$. Multiply by $M^{2k}$, a bound for the number of pairs of words, and choose an integer $d=d_{\rm sep}>X+2\log M/\theta_P$. The exponent is $b_P=\theta_P(d_{\rm sep}-X)-2\log M>0$. Enlarging the constant gives the same bound for the remaining finitely many values of $k$.
\end{proof}

\subsection{Extraction of a convolution}
If the stationary slope-averaged entropy rate falls below the symbolic value $Q$, exponentially separated cylinder translations retain entropy below the cylinder scale. To extract this entropy as a convolution factor, we view the slope and translation of a cylinder as the coefficients of an affine function of time. Evaluation at one time in each of two separated intervals recovers both coefficients with a uniform inverse bound. The following lemma makes this step quantitative.
\begin{lemma}[A positive-entropy factor with negligible convolution growth]\label{extraction}
Fix $P$ of rank $p\ge1$ and suppose $\kappa=\kappa(P)<Q$. There exist $c>0$, an integer $n_1$ and, for every $n\ge n_1$, a slope $a_n\in B$ and a probability measure $\tau_n$ on $[-2/S,2/S]^p$ such that $H_n(\tau_n)\ge cn$ for every $n\ge n_1$ and, as $n\to\infty$,
\begin{equation}\label{extract}
 h_n(\nu_{a_n})\to\kappa,\qquad
 h_n(\tau_n*\nu_{a_n})\to\kappa .
\end{equation}
\end{lemma}
\begin{proof}
\emph{Order of scales.} Fix the observation length $n$. The outer scale $\ell$ will tend to infinity with $n$ fixed. It determines $k=\lceil\ell/\ell_S\rceil$ and $K=d_{\rm sep}k+2$. Sample $i$ uniformly from $0,\ldots,K-1$. Its prefix length $k_i=\lceil i/\ell_S\rceil$ is determined by $i$. In Step~3 we restrict to $i\ge I(n)$ so that replacing the time coordinate by a cylinder endpoint changes it by at most $2^{-n}$. This threshold is independent of $\ell$. Step~4 chooses $\ell(n)$ before taking $n\to\infty$.

\emph{1. Entropy in the coefficients.}
Fix $a$ sampled from $\eta$. A random length-$k$ prefix $w$ determines the pair $(\zeta,\upsilon)=(R_w a,S^{-k}c_w(a))$. Denote its distribution at fixed $a$ by $\Theta_a$. Let $\widehat{\nu}_{a,k}=(z\mapsto S^{-k}z)_*\nu_a$
be the rescaled projected measure at this fixed prefix length $k$. By \eqref{cyl} this is the distribution of $\mathcal Y=PZ'+\zeta T'+\upsilon$ for an independent tail. Given $\mathcal D_0(\mathcal Y)$, the pair $(\zeta,\upsilon)$ lies in at most $J_0=O(1)$ possible unit cells, since $\|\zeta\|\le B_0$ and $\|\mathcal Y-\upsilon\|\le B_Y$. Hence $H_0(\Theta_a)\le H_0(\widehat{\nu}_{a,k})+\log J_0\le H_\ell(\nu_a)+C_{\rm grid}(r,1/S)+\log J_0$ by Lemma~\ref{grid}. For slopes outside the exceptional set in Lemma~\ref{sep}, distinct words occupy distinct level-$K$ cells of $\Theta_a$, since $S^{-k}\ge1$. There $H_K(\Theta_a)=H(W)=k\hw$: the separated atoms have the original product weights. For slopes in the exceptional set, $H_K(\Theta_a)\ge0$. Consequently
\[
 \int (H_K-H_0)(\Theta_a)\,d\eta
 \ge k\hw(1-O_P(2^{-b_Pk}))-\mathcal E_\ell(P)-O(1)
 =(Q-\kappa)\ell+o(\ell)
\]
by the definition of $Q$ and Lemma~\ref{loc}. As $K/\ell\to d_{\rm sep}/\ell_S$, the left side is at least $\beta_{\rm ent} K$ for all large $\ell$, with $\beta_{\rm ent}=(Q-\kappa)\ell_S/(2d_{\rm sep})>0$ independent of $n$.

Choose a level-$i$ coefficient cell according to its mass. Denote its conditional measure before rescaling by $\xi$ and select a representative $(\zeta_0,\upsilon_0)$ from it (use the first word in lexicographic order). Put $\mathfrak a=H_{i+n}(\xi)/n$. The cell has at most $2^{2pn}$ level-$i+n$ children. By \eqref{tele} applied to $\Theta_a$, $0\le \mathfrak a\le2p$ and $\E \mathfrak a\ge\beta_{\rm ent}-2pn/K$. Once $2pn/K\le\beta_{\rm ent}/2$,
\[
 \beta_{\rm ent}/2\le \beta_{\rm ent}/4+2p\,\mathbb P(\mathfrak a\ge\beta_{\rm ent}/4),
 \qquad \mathbb P(\mathfrak a\ge\beta_{\rm ent}/4)\ge p_{\rm coef}:=\beta_{\rm ent}/(8p).
\]

\emph{2. Evaluation at two times.}
Let $i_-,i_+$ be the leftmost and rightmost horizontal branches, and let
$\mathcal O_0,\mathcal O_1$ be the closed intervals of the words $i_-i_-$ and $i_+i_+$.
The interval for $i_-i_+$ lies between them, so their distance is at least $g_t=\rho^2$.
Whenever $k_i\ge2$, an independently sampled prefix starts with each of these two words with probability at least
$p_t=(\min_j\pi_j)^2$. Its left endpoint then belongs to the corresponding interval. Define
\[
 \tau_{\xi,t}=\big((\zeta,\upsilon)\mapsto S^{-k_i}((\zeta-\zeta_0)t+\upsilon-\upsilon_0)\big)_*\xi .
\]
The support of $\tau_{\xi,t}$ lies in $[-2/S,2/S]^p$ for $0\le t\le1$, since the coefficient pair is confined to a cell of side $2^{-i}$ (in each coordinate). For $t\in \mathcal O_0$ and $t'\in \mathcal O_1$ in $[0,1]$, the two evaluations recover the slope by subtraction and division by $t'-t$, and then the intercept by substitution. In maximum norm the two evaluations have forward norm at most $2$ and inverse norm at most $1+2/g_t$: the slope has norm at most $2/g_t$ times the maximum observed norm, and either evaluation recovers the intercept. Rescale the input by $2^i$. The remaining scale factor $2^{-i}S^{-k_i}$ is in $[1,1/S)$. Lemma~\ref{grid}(2) and entropy subadditivity thus give, with a constant $C_{\rm ev}$ independent of $i,\xi,t,t'$,
\[
 H_n(\tau_{\xi,t})+H_n(\tau_{\xi,t'})\ge n\mathfrak a-C_{\rm ev}.
\]
For $\mathfrak a\ge\beta_{\rm ent}/4$ and $C_{\rm ev}<\beta_{\rm ent} n/20$ there cannot be a time of entropy below $\beta_{\rm ent} n/10$ in each interval. Hence one entire interval (intersected with $[0,1]$) gives entropy at least $cn$, $c=\beta_{\rm ent}/10$. The prefix $u$ is independent of the coefficient cell given $i,a$.

\emph{3. Linearization at the entropy mesh.}
Take $I(n)=\max(\lceil2\ell_S\rceil,\lceil n\ell_S/X\rceil)$. Then $i\ge I(n)$ gives $\rho^{k_i}\le2^{-n}$. Write $\mathcal I_n=\{i\ge I(n)\}$. Choose $\ell$ large enough that $I(n)/K\le p_{\rm coef}/2$. Steps 1 and 2 give
\begin{equation}\label{pevent}
 \mathbb P(E_n\mid \mathcal I_n)\ge p_t p_{\rm coef}/2=:p_*>0,\qquad
 E_n=\{H_n(\tau_{\xi,x_u})\ge cn\}.
\end{equation}
In fact $\mathbb P(\mathfrak a\ge\beta_{\rm ent}/4,\mathcal I_n)\ge p_{\rm coef}-I(n)/K$, and each such choice has probability at least $p_t$ of the required $u$.

Condition on the coefficient cell and the prefix $u$ of the independent tail. Split $\mathcal Y$ as
\[
 [(\zeta-\zeta_0)T'+\upsilon-\upsilon_0]+[PZ'+\zeta_0T'+\upsilon_0].
\]
Replace $T'$ by $x_u$ in the first bracket. On $\mathcal I_n$ the difference in maximum norm is at most $2^{-i}\rho^{k_i}\le2^{-(i+n)}$. By Lemma~\ref{grid}(1) the conditional entropy at level $i+n$ changes by at most $p\log5$ bits. After replacement the brackets are conditionally independent: one uses only the coefficient pair in the chosen cell, the other only the independent tail with prefix $u$. By \eqref{cyl} their sum is a translate of the dilation by $S^{k_i}$ of $\tau_{\xi,x_u}*\nu_{R_u(\zeta_0)}$. Lemma~\ref{grid}(2) compares its fine entropy with $H_n$ of this convolution at cost $C_{\rm grid}(r,1/S)$. The support of the original conditional projected measure has diameter at most $C_*2^{-i}$: the first bracket has maximum norm at most $2\cdot2^{-i}$, and the second has diameter at most $2B_Y S^{k_i}$. Thus a further cost bounded by $p\log(\lceil C_*\rceil+2)$ covers its coarse index. Taking conditional expectations in \eqref{mix} gives on each retained scale, with a fixed $C_{\rm lin}$ equal to the sum of these three bounds,
\begin{equation}\label{lin}
 H_{i+n}(\widehat{\nu}_{a,k})-H_i(\widehat{\nu}_{a,k})\ge
 \E_{\xi,u} H_n(\tau_{\xi,x_u}*\nu_{R_u(\zeta_0)})-C_{\rm lin}.
\end{equation}
Here $C_*$ may be taken as $4+2B_Y$ in maximum norm. All three error bounds are uniform in $n$.

\emph{4. Simultaneous selection.}
Let $\mathfrak f_n=h_n(\tau_{\xi,x_u}*\nu_{R_u(\zeta_0)})$ and $\mathfrak d_{i,n}=(H_{i+n}-H_i)(\widehat{\nu}_{a,k})/n$. Averaging $\mathfrak d_{i,n}$ first over all $i$ and then over $a$ gives
\[
 \E\mathfrak d_{i,n}=\frac{\mathcal E_{\ell+K}(P)-\mathcal E_\ell(P)}K+O(pn/K+C_{\rm grid}(r,1/S)/K).
\]
This identity follows by telescoping with \eqref{tele} and then comparing the two endpoint entropies of $\widehat{\nu}_{a,k}$ with the shifted grids of $\nu_a$. Since $\mathcal E_j(P)=\kappa j+o(j)$ and $K/\ell$ tends to a positive constant, the right side tends to $\kappa$ for fixed $n$. Since $\mathfrak d_{i,n}\ge0$, \eqref{lin} implies $\E(\mathfrak f_n\mid \mathcal I_n)\le\E\mathfrak d_{i,n}/\mathbb P(\mathcal I_n)+C_{\rm lin}/n$. For each sufficiently large $n$ now take $\ell(n)$ large enough for Steps 1--3, $I(n)/K\le1/n$ and $|\E\mathfrak d_{i,n}-\kappa|\le1/n$. These choices of $\ell(n)$ are made with $n$ fixed.

Write $C=C_{\rm grid}(r,1)$ and set
\[
 e_n=\max\{0,\kappa-\inf_{a\in B}h_n(\nu_a)\}+C/n .
\]
Lemma~\ref{loc} gives $e_n\to0$. By Lemma~\ref{grid}(3), $\mathfrak f_n\ge\kappa-e_n$ for every outcome. The choice of $\ell(n)$ gives $\E(\mathfrak f_n\mid\mathcal I_n)\le\kappa+o(1)$, and \eqref{pevent} gives $\mathbb P(E_n\mid\mathcal I_n)\ge p_*>0$. Hence
\[
 0\le\E(\mathfrak f_n-\kappa+e_n\mid\mathcal I_n\cap E_n)
 \le p_*^{-1}\E(\mathfrak f_n-\kappa+e_n\mid\mathcal I_n)=o(1).
\]
Choose an outcome in $\mathcal I_n\cap E_n$ with excess at most its conditional mean plus $1/n$. For $\tau_n=\tau_{\xi,x_u}$ and $a_n=R_u(\zeta_0)$ this gives $H_n(\tau_n)\ge cn$ and $\mathfrak f_n\to\kappa$. Since $\kappa-e_n\le h_n(\nu_{a_n})\le\mathfrak f_n+C/n$ by Lemma~\ref{grid}(3) and the definition of $e_n$, both limits in \eqref{extract} hold for the same pair.
\end{proof}

\subsection{The inverse theorem and induction on the target dimension}\label{inverse_sec}
For a subspace $V\le\R^p$ let $\operatorname{pr}_{V^\perp}$ denote orthogonal projection, with image viewed inside $\R^p$. A probability measure $\theta$ is $(V,\varepsilon)$-concentrated if some translate of $V$ has open $\varepsilon$-neighborhood of mass at least $1-\varepsilon$. It is $(V,\varepsilon,m)$-saturated if
\[
 h_m(\theta)\ge\dim V+h_m((\operatorname{pr}_{V^\perp})_*\theta)-\varepsilon .
\]
Here the projected measure is viewed in $\R^p$ with the ambient dyadic grid. In orthonormal coordinates on $V^\perp$ its normalized entropy differs by at most $C_p/m$, uniformly in $V$ and in the measure~\cite[Section~3.3]{hoch}. Indeed, on $V^\perp$ each cell of either grid has diameter at most $\sqrt p\,2^{-m}$ and meets at most $(\lceil\sqrt p\rceil+2)^p$ cells of the other grid. The partition entropy comparison gives $C_p=p\log_2(\lceil\sqrt p\rceil+2)$.

\begin{littheorem}[Large-block form of Hochman~{\cite[Theorems 2.8 and 4.17]{hoch}}]\label{litH}
Fix an integer $p\ge1$, $R>0$ and $\varepsilon>0$. There is an integer $m_0=m_0(p,R,\varepsilon)\ge1$ such that for every integer $m\ge m_0$ there exist $\delta>0$, depending on $p,\varepsilon,R,m$, and then a threshold $n_0$, depending on $p,\varepsilon,R,\delta,m$, such that for every $n\ge n_0$ the following holds for every pair of probability measures $\omega,\tau$ on $[-R,R]^p$. If $h_n(\omega*\tau)<h_n(\omega)+\delta$, there are subspaces $V_0,\ldots,V_n$ such that the probability of both events
\[
 \omega^{x,i}\text{ is }(V_i,\varepsilon,m)\text{-saturated},\qquad
 \tau^{y,i}\text{ is }(V_i,\varepsilon)\text{-concentrated}
\]
exceeds $1-\varepsilon$. Here $i$ is uniform on $\{0,\ldots,n\}$, $x\sim\omega$ and $y\sim\tau$ independently of each other and of $i$. The components are conditioned on their ambient dyadic cells and mapped back to the unit cube, as in \eqref{inc}.
\end{littheorem}
Hochman's Theorem~2.8 uses the concentration, saturation and component conventions of Definitions~2.1, 2.3 and~2.4 in~\cite{hoch}. Theorem~4.17 in Section~4.6 of the same paper restates the result for measures on $[-R,R]^p$. By Proposition~3.19, Hochman reduces the proof to the case of sufficiently large $m$. This is the case used below. Passing between intrinsic and ambient grids costs at most $C_p/m$ in normalized entropy. Apply the theorem at tolerance $\varepsilon/2$ and choose $m_0$ with $C_p/m_0\le\varepsilon/2$. The resulting saturation and concentration events then have joint probability greater than $1-\varepsilon$, uniformly over the measures in the stated box.

\begin{proposition}[Entropy of the projections]\label{kappaprop}
For every $0\le p\le r$ and every $P:\R^r\to\R^p$ with orthonormal rows, $\kappa(P)=\min(p,Q)$.
\end{proposition}
\begin{proof}
The proof is by induction on $p$. Rank zero has $\kappa=0$. Suppose the assertion holds at all lower ranks, fix $P$ of rank $p\ge1$, write $Q_e=\min(e,Q)$ for $0\le e\le p$ and assume for a contradiction that $\Delta_{\rm ent}=Q_p-\kappa(P)>0$. Lemma~\ref{extraction} applies because $\kappa(P)<Q$. Write $\omega_n=\nu_{a_n}$. Its support and that of $\tau_n$ lie in $[-R,R]^p$ with $R=1+\max(B_Y,2/S)$. We shall obtain a gain of $\Delta_{\rm ent}$, up to errors, on a positive fraction of the scales, while retaining the lower bound $\kappa(P)-\gamma$ on the other late scales. The tolerances and block length are chosen at the end of the proof. For each admissible $\varepsilon,m$, Theorem~\ref{litH} supplies $\delta,n_0$ uniformly over measures in this box. By \eqref{extract} its strict entropy hypothesis holds for all sufficiently large $n\ge n_0$. For each such $n$, apply the inverse theorem once to $(\omega_n,\tau_n)$ and use its subspaces $V_i$ for both measures at scale $i$.

\emph{Proportion of scales with a nonzero saturation space.}
Put $K=n+1$ and $u=\sqrt{\varepsilon}$ with $\varepsilon<1/8$. Average over all $K$ scales $0,\dots,n$. Call a scale good if, conditional on that scale, the probability that either conclusion of Theorem~\ref{litH} fails is at most $u$. By Markov's inequality, the fraction of other scales is at most $u$. At a good scale, concentration and saturation therefore each hold with probability at least $1-u$.

Let $\ell=\lfloor\tfrac12\log(1/\varepsilon)\rfloor$. This auxiliary length $\ell$ is used only in the entropy increment estimate for the convolution factor $\tau_n$. The inverse theorem is invoked at the main length $m$. On a good scale with $V_i=0$, a concentrated component of $\tau_n$ assigns at least $1-\varepsilon$ of its mass to a ball of radius $\varepsilon$. This ball meets at most $4^p$ level-$\ell$ cubes since $2^{-\ell}\ge\sqrt\varepsilon$. The indicator of membership in the ball has entropy at most one bit, and the entire component is contained in $[0,1)^p$ with $2^{p\ell}$ cells. Its entropy is at most $1+2p+\varepsilon p \ell$. The exceptional components have entropy at most $p\ell$. The mean normalized entropy on this scale is consequently at most $a_\varepsilon=(1+2p)/\ell+p(\varepsilon+u)$. Every other scale has mean at most $p$. Write $f_n$ for the fraction of good scales with $V_i\ne0$, and $C_R=p\log(\lceil2R\rceil+2)$ for a unit-grid entropy bound. Lemma~\ref{increments} at length $\ell$ gives
\[
 c-\frac{C_R+p\ell+c}{K}\le a_\varepsilon+p(f_n+u)
\]
because $H_n(\tau_n)\ge cn$. Whenever $a_\varepsilon+pu\le c/4$ and $(C_R+p\ell+c)/K\le c/4$, this gives
\begin{equation}\label{proportion}
 f_n\ge f_*:=c/(2p)>0 .
\end{equation}
This lower bound is independent of $m$ and of all subsequent decreases of $\varepsilon$.

\emph{Entropy gain at such a scale.}
Let $i_0(m)$ be the threshold in Lemma~\ref{loc} at tolerance $\gamma$. By induction $\kappa_e(P)=Q_e$ for every $e<p$. On a good scale $i\ge i_0(m)$ with $v=\dim V_i>0$, choose any orthonormal coordinates $U_i:\R^p\to\R^{e}$ on the orthogonal complement ($e=p-v$, kernel $V_i$). On the saturated components,
\[
 h_m(\omega_n^{x,i})\ge v+h_m((U_i)_*\omega_n^{x,i})-\varepsilon-C_{\rm grid}(r,1)/m .
\]
The projected image of a unit cube has Euclidean diameter at most $\sqrt p$, so its level-$m$ entropy is at most $em+c_r$, with $c_r$ depending only on $r$. Dropping the nonsaturated components, whose proportion is at most $u$, loses at most $u(p+c_r/m)$ from the mean of $v+h_m((U_i)_*\omega_n^{x,i})$. Lemma~\ref{frames} applies to the full mean of this projected entropy, at this particular $U_i$ and $a=a_n$. Thus
\begin{align}
 \E_x h_m(\omega_n^{x,i})&\ge v+Q_e-\varepsilon_{\rm gain}\ge Q_p-\varepsilon_{\rm gain},\label{gain}\\
 \varepsilon_{\rm gain}&=\gamma+\varepsilon+C_{\rm grid}(r,1)/m+u(p+c_r/m).\notag
\end{align}
The four terms bound the errors from the lower-rank estimate, saturation, the change of grid, and the discarded components, respectively. The inequality $v+\min(p-v,Q)\ge\min(p,Q)$ covers also $v=p$. Every other scale $i\ge i_0(m)$ has the lower bound $\E_x h_m(\omega_n^{x,i})\ge\kappa-\gamma$ by Lemma~\ref{loc}, whether or not the scale is good.

At most $i_0(m)$ scales have $i<i_0(m)$, and their entropies are nonnegative. Let $f'_n$ be the fraction of all $K$ scales that are good and satisfy $i\ge i_0(m)$ and $V_i\ne0$. By \eqref{proportion}, $f'_n\ge f_*-i_0(m)/K$. Averaging the last two bounds gives
\[
 \frac1K\sum_{i<K}\E_x h_m(\omega_n^{x,i})
 \ge \kappa-\gamma+ f'_n\Delta_{\rm ent}-\varepsilon_{\rm gain}-p\,i_0(m)/K ,
\]
where the sum starts at zero. By \eqref{tele} the difference of the left side from $h_n(\omega_n)$ is at most $(pm+C_R+2p)/K$. Here $0\le H_{n+1}-H_n\le p$ and $H_n\le H_0+pn$ for the projected measure. At fixed $m$, \eqref{extract} and passage to the limit therefore yield $f_*\Delta_{\rm ent}\le\gamma+\varepsilon_{\rm gain}$. Choose the parameters in the following order. First take $2\gamma<f_*\Delta_{\rm ent}/4$, take $0<\varepsilon<1/8$ small enough that $a_\varepsilon+pu\le c/4$ and $\varepsilon+up<f_*\Delta_{\rm ent}/4$, then enlarge $m$ until $(C_{\rm grid}(r,1)+uc_r)/m<f_*\Delta_{\rm ent}/4$, Lemma~\ref{loc} applies both with $U=I_p$ and with $U=U_i$, and $m\ge m_0(p,R,\varepsilon)$. With $\varepsilon$ and $m$ fixed, take $\delta,n_0$ from Theorem~\ref{litH} and let $n\to\infty$. These choices contradict the preceding inequality, proving $\kappa(P)=Q_p$ at rank $p$.
\end{proof}

\begin{proof}[Proof of Theorem~\ref{hom}(1)]
The slope radius $B_0$ fixed before Lemma~\ref{loc} is at least the radius in the statement, so the slopes in the statement lie in $B$. By Proposition~\ref{kappaprop} and the last assertion of Lemma~\ref{loc}, $\liminf_{n\to\infty}\inf_{a\in B}h_n(\nu_a)\ge\min(p,Q)$. The cylinder count in the proof of Lemma~\ref{loc} gives $H_i(\nu_a)\le i\min(p,Q)+O(1)$ uniformly for $a\in B$. The two bounds give the asserted uniform convergence.
\end{proof}

\begin{lemma}[All slopes and all target directions]\label{verticalsets}
Let $\Lambda$ be the attractor of the homogeneous system in
Theorem~\ref{hom}, and put
\[
 K_{P,a}=\{Pz+at:(t,z)\in\Lambda\}.
\]
For every $1\le p\le r$, every orthonormal-row map
$P:\R^r\to\R^p$, and every $a\in\R^p$,
\[
 \dim_{\mathrm H}K_{P,a}
 =\underline{\dim}_{\B}K_{P,a}
 =\overline{\dim}_{\B}K_{P,a}
 =\min\left(p,\frac{\log M}{\log(1/S)}\right).
\]
\end{lemma}
\begin{proof}
With uniform weights put $s=\min(p,\log_2 M/\ell_S)$ and choose an invariant slope ball containing $a$. Theorem~\ref{hom}(1) gives at least $2^{(s-\varepsilon)m}$ occupied dyadic cells at $m=\lfloor k\ell_S\rfloor$, uniformly over slopes in that ball for large $k$. Select one support point in each cell. Keep a largest class of cell indices modulo an integer $D>4+4B_Y$ in each coordinate, where $B_Y$ bounds all projected supports in maximum norm. Chosen points are at least $(D-1)2^{-m}$ apart. By \eqref{cyl} each belongs to a length-$k$ cylinder of diameter at most $2B_Y S^k$. The chosen cylinders have mutual distance at least $S^k$.
Repeat at each descendant slope with exactly $N_k=\lfloor D^{-p}2^{(s-\varepsilon)m}\rfloor$ children per parent. At depth $j$ different pieces have separation at least $S^{kj}$. By packing, a ball of radius $S^{kj}$ meets at most a fixed number of these pieces. Equal masses thus give Hausdorff dimension at least $\log_2 N_k/(k\ell_S)$, by comparison of successive radii. Let $k\to\infty$ and $\varepsilon\downarrow0$. On the other hand, $M^k$ sets of diameter $O(S^k)$ cover the projection, and the ambient dimension is $p$.
\end{proof}

\section{Different branch ratios via permutation blocks}\label{permsec}
For a strictly positive probability vector $\mathbf p$ on $J\ne\varnothing$ put $\hw_{\bf p}=-\sum p_i\log p_i$, $\chi_t=-\sum p_i\log b_i$ and $\chi_v=-\sum p_i\log a_i$. Here $\hw_{\bf p}$ is the entropy of the Bernoulli measure with weights $\mathbf p$ for the shift, and $\chi_t,\chi_v$ are the horizontal and vertical Lyapunov exponents, with the positive sign convention for contractions. All three use natural logarithms. In Lemma~\ref{homgraph} the logarithms are to base two.

\begin{lemma}[Two-scale lower bound for a homogeneous graph]\label{homgraph}
Under the hypotheses of Theorem~\ref{hom}, let $\Lambda$ be the attractor and put $Q_0=(\log_2 M)/\ell_S$. Then
\[
 \dim_{\rm H}\Lambda\ge(\log_2 M)/X+(1-\ell_S/X)\min(r,Q_0).
\]
\end{lemma}
\begin{proof}
Use base-two logarithms as in Section~\ref{sec:homogeneous}. If $Q_0\le r$, Lemma~\ref{verticalsets} with $P=I_r$ and slope zero gives a Lipschitz image of dimension $Q_0$. The right-hand side equals $Q_0$.

Suppose $Q_0>r$, use the maximum norm and take $0<\varepsilon<1$. Use uniform symbol weights,
$P=I_r$ and an invariant slope ball containing zero. Fix a large block length $k$ and put
$m=\lfloor k\ell_S\rfloor$. For a uniform word $w$ of length $k$ the translation $c_w(a)$
is within $B_Y S^k$ of a coupled random point with distribution $\nu_a$. Thus their dyadic entropies at
level $m$ differ by $O(1)$, uniformly in the slope $a$ and in $k$, by bounded grid
multiplicity. The translations occupy at most $C2^{rm}$ cells. By Theorem~\ref{hom}(1)
their entropy is at least $rm-o(m)$ uniformly in $a$. If cells of mass greater than
$2^{-m(r-\varepsilon)}$ have total mass $v$, the binary split into these cells and the others
bounds the entropy above by $1+v(r-\varepsilon)m+(1-v)(rm+\log_2 C)$.
Consequently $v\le1/2$ for all sufficiently large $k$.

At each slope retain $n_k=\lceil M^k/2\rceil$ words whose translations belong to the
remaining cells. A ball of radius $S^k$ meets at most
$C_1M^kS^{k(r-\varepsilon)}$ retained projected cylinders. Indeed, their translations lie in
a ball of radius $(1+B_Y)S^k$, which meets boundedly many cells, each containing at most
$M^k2^{-m(r-\varepsilon)}$ retained words. The constant $C_1$ is independent of $k,a,\varepsilon$.
Iterate the choice at the slope $R_w0$ of each retained prefix $w$. Give each of its
chosen children conditional probability $1/n_k$ and push the resulting measure forward to $\Lambda$ under the coding map.

For a ball of radius $0<u<1$ let $j,j'$ be the unique nonnegative integers with
$\rho^{k(j+1)}<u\le\rho^{kj}$ and $S^{k(j'+1)}<u\le S^{kj'}$. Then $j\le j'$.
At most four retained level-$j$ time intervals meet the ball, including endpoint
contacts: their interiors are disjoint and their lengths are $\rho^{kj}\ge u$.
At each subsequent depth $i=j,\ldots,j'-1$, the normalized vertical radius is
$u/S^{ki}\le S^k$. The preceding count applies inside each parent by \eqref{cyl}.
The ball therefore has mass at most
\[
 4(2/M^k)^j\big(2C_1 S^{k(r-\varepsilon)}\big)^{j'-j}.
\]
For fixed $\varepsilon,k$, rounding the scale indices $j,j'$ introduces a multiplicative constant
independent of $u$. The mass distribution principle gives dimension at least
\[
 \frac{\log_2 M-1/k}{X}
 +\left(1-\frac{\ell_S}{X}\right)
 \left(r-\varepsilon-\frac{\log_2(2C_1)}{k\ell_S}\right)
\]
for all sufficiently large $k$. Each constructed measure is supported on $\Lambda$;
letting $k\to\infty$ and then $\varepsilon\downarrow0$ proves the claim.
\end{proof}

\begin{lemma}[Permutation blocks preserve the shear span]\label{perms}
Let $r\ge1$. Assume \eqref{rank} and $A(1)>1$, so that at least two branches are active. For every strictly positive probability vector $\mathbf p$ on $J$ with $\chi_t>\chi_v$ there is a sequence of homogeneous subsystems of the normalized graph satisfying the hypotheses of Theorem~\ref{hom}. Their attractors give
\begin{equation}\label{varlower}
 \HH G\ge \hw_{\bf p}/\chi_t+(1-\chi_v/\chi_t)\min(r,\hw_{\bf p}/\chi_v).
\end{equation}
The right-hand side of \eqref{varlower} is the Lyapunov dimension of the self-affine system on $\R\times E$ with Bernoulli weights $\mathbf p$. In increasing order, the positive contraction exponents are $\chi_v$ with multiplicity $r$, followed by $\chi_t$.
\end{lemma}
\begin{proof}
For an integer $\ell\ge1$ choose $m_i=2\lceil \ell p_i\rceil$, put $k=\sum_i m_i$ and form the set of all $M_k=k!/\prod_i m_i!$ distinct permutations of this multiset. Every block sends the normalized graph into itself. Concatenating infinitely many blocks gives a code in the original system, so their attractor lies in that graph. Different blocks have disjoint horizontal interiors by the first index where they differ. Their horizontal ratio is $\rho_k=\prod b_i^{m_i}$. Their signed vertical product is positive since every $m_i$ is even, and equals $S_k=\prod a_i^{m_i}$. As $\ell\to\infty$ one has $m_i/k\to p_i$. Stirling's formula $\log n!=n\log n-n+O(\log(n+1))$ therefore gives
\[
 (k^{-1}\log M_k,-k^{-1}\log\rho_k,-k^{-1}\log S_k)\longrightarrow
 (\hw_{\bf p},\chi_t,\chi_v).
\]
In particular $0<\rho_k<S_k<1$ eventually.

By \eqref{pieceeq}, the shear of a block $w$ is $\Delta_w$, the coefficient of $t$. Fix distinct $i,j\in J$ and form a prefix $u$ from the multiset after removing one copy of each of $i,j$. Both $uij$ and $uji$ belong to the block alphabet. Applying \eqref{piecerec} twice gives
\[
 \Delta_{uij}-\Delta_{uji}
 =s_u\big[(b_j-s_j)\Delta_i-(b_i-s_i)\Delta_j\big].
\]
Every symbol of $u$ is active, so $s_u\ne0$. As $i,j$ vary, \eqref{rank} shows that the differences of block shears span $E$. Apply Lemma~\ref{homgraph} to each block attractor, which lies in the normalized graph. The limits above give the Bernoulli entropy and the two Lyapunov exponents. Their ratios are independent of the logarithm base. Taking the limit gives \eqref{varlower}. The constants may depend on the fixed block system.
\end{proof}

\begin{lemma}[Relative entropy]\label{relent}
Let $B_J=\sum_{i\in J}b_i\le1$. For every strictly positive probability vector $\mathbf p$ on $J$,
\[
 \chi_t-\hw_{\bf p}\ge-\log B_J\ge0 .
\]
\end{lemma}
\begin{proof}
This is nonnegativity of the relative entropy of $\mathbf p$ with respect to the probability vector $(b_i/B_J)_{i\in J}$, namely $\sum_{i\in J}p_i\log(p_iB_J/b_i)\ge0$.
\end{proof}

\begin{lemma}[Dimension of the graph]\label{graphvalue}
For interpolation data of positive residual rank $r$, if $A(1)\le1$ or \eqref{rank} holds, the Hausdorff and both Minkowski dimensions of $G$ equal $D_r$.
\end{lemma}
\begin{proof}
The projection of $G$ onto the time axis is $[0,1]$, so $\HH G\ge1$. With the graph upper bound of Section~\ref{upperbound} this gives $1\le\HH G\le\underline{\dim}_{\B}G\le\overline{\dim}_{\B}G\le D$. This settles the case $A(1)\le1$, where $D=1$. Let $A(1)>1$. The rank hypothesis now gives Lemma~\ref{perms}. In the middle regime $A(r)\le1$ take $p_i=a_i^q$, which is a probability vector because $A(q)=1$. Then $\hw_{\bf p}=q\chi_v>\chi_v$, hence $\chi_t>\chi_v$ by Lemma~\ref{relent}, and since $q\le r$ the bound \eqref{varlower} equals $q$. In the last regime $A(r)>1$ take $p_i=a_i^r b_i^\delta$. Then $\hw_{\bf p}=r\chi_v+\delta \chi_t$, so Lemma~\ref{relent} gives $\chi_t\ge r\chi_v/(1-\delta)>\chi_v$ and $\hw_{\bf p}/\chi_v>r$. The bound \eqref{varlower} equals $r+\delta$. In both regimes $\HH G\ge D$.
\end{proof}

\section{Proof for arbitrary linear images}\label{rangesection}
\begin{proof}[Proof of Theorem~\ref{main}]
For positive residual rank, Section~\ref{upperbound} and Lemma~\ref{graphvalue} give the upper bound and the equality in part~\textup{(a)}, respectively. At rank zero the graph is a line segment and $D_0=1$. We prove part~\textup{(b)} for an arbitrary linear map.
Write $\Pi(t,y)=tv+Ly$, $z_i=\Pi(t_i,y_i)$ and $z_*=z_N-z_0=v+L(y_N-y_0)$. Subtracting the endpoint line of $f$ gives
\[
 F(t):=\Pi(t,f(t))=z_0+t z_*+H(t),\qquad H=Lh .
\]
Applying $L$ to the normalized branch equation for $h$ gives the same equation for $H$, with residual interpolation values $Lu_i$ and increments $L\Delta_i$. Thus $F$ is the interpolant of the $z_i$ with the original widths and scaling factors. Its residual space is $V=L(E)$ and its rank generators span $L(W)$, where $W$ is the span in \eqref{rank}.

Since $z_i-z_0=Lu_i+t_i z_*$ and $z_*=z_N-z_0$ occurs among these differences, their span is $V+\R z_*$. Hence $j=\dim V$, while $g=j$ if $z_*\in V$ and $g=j+1$ otherwise. The image $\Pi(G)$ is the range of $F$.

If $j=0$ the image is the segment from $z_0$ to $z_N$, including a singleton when $z_*=0$. Since $D_0=1$ and $g\in\{0,1\}$, its dimension is $g=\min\{D_0,g\}$. For $j\ge1$, the range of $F$ is a Lipschitz image of its graph and lies in $z_0+V+\R z_*$. The graph bound proved in Section~\ref{upperbound}, applied to the transformed data, and the affine dimension bound give the upper bound in part~\textup{(b)}.

Assume now that $A(1)\le1$ or $L(W)=V$. In the latter case the transformed data satisfy \eqref{rank}. If $z_*\notin V$, the map $(t,z)\mapsto z+t z_*$ is a linear isomorphism from $\R\times V$ onto $V+\R z_*$. Up to translation, this map is a bi-Lipschitz equivalence between the graph of $H$ and the image. Lemma~\ref{graphvalue}, applied to the normalized transformed data, gives $D_j$, and $D_j\le j+1=g$.

Suppose $z_*\in V$ and $j\ge1$. If $A(1)\le1$, the image is connected and contains distinct projected interpolation points because $g=j\ge1$. Its Hausdorff dimension is at least one (project onto the line through two such points). Now let $A(1)>1$, and let $A(q)=1$. In the transformed system on $V$ use the frequencies $p_i=a_i^q$. Their entropy is $q\chi_v$, and Lemma~\ref{relent} gives $\chi_t\ge q\chi_v>\chi_v$. Lemma~\ref{perms} constructs homogeneous subsystems on the normalized graph of $H$ whose shear differences span $V$. In orthonormal coordinates on $V$, apply Lemma~\ref{verticalsets} with $P=I_j$ and slope $z_*$. The translated projected attractors lie in $F([0,1])$ and have Hausdorff dimensions $\min(j,\log M_k/(-\log S_k))$, tending to $\min(j,q)$. By \eqref{Ddef} this is $\min(j,D_j)$. Together with the upper bound, this proves \eqref{allrangeformula}.
\end{proof}

\paragraph{Ranges and coordinate graphs.}
Taking $\Pi(t,y)=P_\Gamma y$ gives the range projected onto $\Gamma$. Its residual rank is $\dim P_\Gamma E$ and its affine rank is the affine dimension of the projected values $P_\Gamma y_i$. Taking $\Gamma=\R^n$ gives the full range. For $U\subset E$ of dimension $1\le p\le r$, take first $\Pi(t,y)=P_Uy$ and then $\Pi(t,y)=P_U(y-t(y_N-y_0))$. The second image is a translate of $P_Uh([0,1])$. In each case $j=g=p$. Under the equality hypotheses,
\begin{equation}\label{allinsidespan}
 \dim_{\rm H}P_Uf([0,1])=\dim_{\B}P_Uf([0,1])
 =\dim_{\rm H}P_Uh([0,1])=\dim_{\B}P_Uh([0,1])=\varrho_p ,
\end{equation}
where $\varrho_p=\min(p,D_p)$ equals $1$ for $A(1)\le1$ and $\min(p,q)$ otherwise. If $E=\R^n$, this applies to every $p$-dimensional subspace. For $\Pi(t,y)=(t,Ly)$ the time component makes the endpoint displacement transverse to $\{0\}\times L(E)$, giving the coordinate-graph formula of Section~\ref{statement}.

\section{Exceptional vertical factors and dimension gaps}\label{sec:exceptional}
\subsection{The span of the active shears}
Fix the interpolation points. Put $\sigma_i=b_i-s_i$ and write $r_W=\dim W$ for the dimension of the span on the left of \eqref{rank}.
On $\R^J$ define $Be_i=\Delta_i$ and $\sigma(x)=\sum_{i\in J}\sigma_i x_i$.

\begin{proposition}[Rank on the active branches]\label{prop:active-test}
If $\sigma=0$, then $W=0$. Otherwise
\[
 W=B(\ker\sigma),\qquad r_W=\operatorname{rank}\binom{B}{\sigma}-1 .
\]
For $r\ge1$, condition~\eqref{rank} holds precisely when the active increments span $E$ and no $\vartheta\in E^*$ satisfies $\sigma_i=\vartheta(\Delta_i)$ for every $i\in J$.
\end{proposition}
\begin{proof}
The generators are $B(\sigma_j e_i-\sigma_i e_j)$. If $\sigma_j\ne0$ for some $j$, these vectors before applying $B$ span $\ker\sigma$, because for $x\in\ker\sigma$ their linear combination with coefficients $x_i$ is $\sigma_jx$.
Projection of the image of $x\mapsto(Bx,\sigma(x))$ onto its last coordinate is surjective with kernel $W\times\{0\}$.
Hence $r_W=r$ precisely when $B$ is onto and $\sigma$ is independent of its rows, equivalently when $\sigma\ne\vartheta B$ for every $\vartheta\in E^*$.
If $\sigma=0$, every generator vanishes and $\vartheta=0$ gives a factorization.
\end{proof}

Suppose all branches are active and $r\ge1$. Then $B$ is onto, and the admissible factors at which \eqref{rank} fails lie exactly in the affine subspace
\[
 \Xi=b+\{(\vartheta(\Delta_i))_{i=0}^{N-1}:\vartheta\in E^*\}
 \ \subset\ \{s:\textstyle\sum_i s_i=1\}.
\]
Indeed, $\sigma=b-s$, and either sign of a functional is allowed.
The increments span $E$ and sum to zero. Thus the parametrization is injective, $\dim\Xi=r\le N-1$, and $\Xi$ contains $b$.
For these fixed interpolation data Theorem~\ref{main} therefore applies at every all-active vector outside $\Xi$, even if $\sum_i s_i=1$. When $r<N-1$, $\Xi$ is a proper affine subspace of this hyperplane. For fixed widths and interpolation values, $\Xi$ describes exactly where the active branches lose rank, under the assumptions above. Outside $\mathcal X_N$, Corollary~\ref{cor:generic} gives a simultaneous guarantee for every choice of interpolation data.

\begin{proposition}[Graph dimension at reduced rank]\label{prop:reduced-graph}
For arbitrary vertical factors and widths,
\[
 D_{r_W}\le\HH G\le\underline{\dim}_{\B}G\le\overline{\dim}_{\B}G\le D_r.
\]
\end{proposition}
\begin{proof}
For $r_W=0$ the time projection gives $D_0=1$. For $r_W>0$, project the
values of $h$ onto $W$. The projected increments span $W$, since the
original increments span $E$. Each generator of \eqref{rank} is fixed
by the projection, so the projected system satisfies \eqref{rank}.
Theorem~\ref{main}(a) gives its graph dimension $D_{r_W}$, with the same
widths and vertical factors. The projection is Lipschitz and the endpoint shear
is bi-Lipschitz. The upper bound also follows from Theorem~\ref{main}(a).
\end{proof}

\subsection{Structure on the exceptional subspace}

\begin{proposition}[The exceptional vectors of factors]\label{excstructure}
Let all branches be active, $r\ge1$, and let $s_i=b_i+\vartheta(\Delta_i)$ for all $i$, where
$\vartheta\in E^*$ and $\vartheta\ne0$.
\begin{enumerate}
\item The functional $\vartheta$ is unique and the span in \eqref{rank} equals $\ker\vartheta$.
\item The function $\psi(t)=t+\vartheta(h(t))$ satisfies
\[
 \psi(t_i+b_it)=\psi_i+s_i\psi(t),\qquad \psi_i=\sum_{k<i}s_k,
\]
with $\psi(0)=0$ and $\psi(1)=1$. Its graph is the attractor of the diagonal maps
$(t,x)\mapsto(t_i+b_it,\psi_i+s_ix)$, and composing $\psi$ with the horizontal coding map gives the coding map of the
self-similar system $\mathcal F_\psi=\{x\mapsto\psi_i+s_ix\}_{0\le i<N}$ on the line.
\item Choose $e\in E$ with $\vartheta(e)=1$ and put $\varpi(z)=z-\vartheta(z)e$ and $\bar h=\varpi\circ h$.
Then $h=(\psi-t)e+\bar h$, so $(t,h(t))\mapsto(t,\psi(t),\bar h(t))$ is the restriction of the linear isomorphism $(t,z)\mapsto(t,t+\vartheta(z),\varpi z)$ from $\R\times E$ onto $\R^2\times\ker\vartheta$. The function $\bar h$ is the interpolant of the interpolation points $(t_i,\varpi u_i)$ with
the factors $s_i$. Its residual space is $\ker\vartheta$, and for $r\ge2$ these data satisfy
\eqref{rank}.
\item For a linear map $\Pi(t,y)=tv+Ly$ the transformed data of Section~\ref{rangesection}
satisfy \eqref{rank} if and only if $\ker(L|_E)\not\subset\ker\vartheta$.
\end{enumerate}
\end{proposition}
\begin{proof}
(1) The increments span $E$, so $\vartheta$ is unique. Since $\sigma=-\vartheta B\ne0$ and $B$ is onto, Proposition~\ref{prop:active-test} gives $W=B(\ker(\vartheta B))=\ker\vartheta$.

(2) Apply $\vartheta$ to $h(t_i+b_it)=u_i+\Delta_it+s_ih(t)$ and add $t_i+b_it$:
\[
 \psi(t_i+b_it)=t_i+\vartheta(u_i)+(b_i+\vartheta(\Delta_i))t+s_i\vartheta(h(t))
 =\psi(t_i)+s_i\psi(t).
\]
Moreover $\psi(t_{i+1})-\psi(t_i)=b_i+\vartheta(\Delta_i)=s_i$ and $\psi(t_0)=0$.

(3) Applying $\varpi$ to the equation of $h$ shows that $\bar h$ satisfies the branch equations
of the interpolation points $(t_i,\varpi u_i)$, whose endpoint values vanish; the interpolant is unique. The new
residual increments are $\varpi\Delta_i$ and span $\varpi(E)=\ker\vartheta$. The generators of \eqref{rank} for the new
data are the images under $\varpi$ of the old ones, so by (1) they span $\varpi(\ker\vartheta)=\ker\vartheta$.

(4) The transformed data have residual increments $L\Delta_i$ and residual space $L(E)$, and their
generators span $L(W)=L(\ker\vartheta)$. Since $\ker\vartheta$ is a hyperplane of $E$, one has
$L(\ker\vartheta)=L(E)$ if and only if $\ker(L|_E)\not\subset\ker\vartheta$.
\end{proof}

\begin{theorem}[Exceptional vertical factors]\label{exceptional}
Suppose all branches are active, $r\ge1$, $A(1)>1$, and \eqref{rank} fails. Write $s=b+(\vartheta(\Delta_i))_i$ as above. Since $A(1)>1$, we have $s\ne b$ and $\vartheta\ne0$. Let $\mathcal F_\psi$ be as in Proposition~\ref{excstructure}.
\begin{enumerate}
\item[\textup{(a)}] For every linear map $\Pi(t,y)=tv+Ly$ with
$\ker(L|_E)\not\subset\ker\vartheta$, formula \eqref{allrangeformula} holds.
\item[\textup{(b)}] If $r\ge2$, then
\[
 D_{r-1}\le\dim_{\rm H}G\le\underline{\dim}_{\B}G\le\overline{\dim}_{\B}G\le D_r .
\]
In particular, if $A(r-1)\le1$, the Hausdorff and both Minkowski dimensions of $G$ equal $D_r$.
\item[\textup{(c)}] If $r=1$ and $\mathcal F_\psi$ is exponentially separated in the sense of~\cite[Definition~1.1]{rap}, the Hausdorff and both Minkowski dimensions of $G$ equal $D_1$.
\end{enumerate}
\end{theorem}
\begin{proof}
(a) By Proposition~\ref{excstructure}(4) the transformed data satisfy the rank condition. The fixed-map criterion following Theorem~\ref{main} gives \eqref{allrangeformula}.

(b) Proposition~\ref{excstructure}(1) gives $r_W=r-1$. Apply Proposition~\ref{prop:reduced-graph}. If $A(r-1)\le1<A(1)$, Definition~\ref{def:D} gives $D_{r-1}=D_r=q$.

(c) Here $E=\R e$, and after subtracting the endpoint line, Proposition~\ref{excstructure}(3) gives a bi-Lipschitz affine identification of $G$ with the graph of $\psi$, the attractor of the planar diagonal system
$(t,x)\mapsto(t_i+b_it,\psi_i+s_ix)$. Its first coordinate system $\{t\mapsto t_i+b_it\}$ is
exponentially separated, because distinct horizontal cylinders of length $m$ have disjoint interiors and lengths at least $(\min_i b_i)^m$. Their left endpoints are therefore at least this far apart. Some $i$ has $a_i\ne b_i$, since $\sum_ia_i>1=\sum_ib_i$.
The theorem of Rapaport on diagonal systems~\cite[Theorem~1.3]{rap} gives
$\dim_{\rm H}G=\min\{2,d_{\rm aff}\}$, where $d_{\rm aff}$ is the affinity dimension. For
$1\le\beta\le2$ the singular value function of the word $w$ equals
$\max\{a_wb_w^{\beta-1},\,b_wa_w^{\beta-1}\}$. Summing over length-$m$ words puts the pressure sum between the larger of the two product sums and twice that value. Taking $m$th roots and passing to the limit gives the maximum of the two one-step sums. Thus $d_{\rm aff}$ is the larger root of
$\sum_ia_ib_i^{\beta-1}=1$ and $\sum_ib_ia_i^{\beta-1}=1$. The second root is $1$ and the first is
$D_1$. The upper bound of Theorem~\ref{main} completes the proof.
\end{proof}

For an instance with $q=r-1$, take four equal widths, interpolation values
$0,e_1,e_2,e_3,0$ in $\R^3$ and vertical factors $(1/2,1/2,1/2,-1/2)$.
Here $r=3$. The only relation among the four increments is their zero
sum, so the row space of $B$ is the hyperplane of zero-sum vectors.
Since $\sum_i(s_i-b_i)=0$, the factors lie on $\Xi$.
Thus $W$ has dimension two, and $4(1/2)^2=1$ gives graph dimension two.

The following criterion uses only the interpolation data and the vertical factors. It is analogous to the uniform fibre condition in the carpets of Bedford and McMullen.

\begin{proposition}[Balanced vertical covering]\label{balanced}
Let $r=1$, let all branches be active, $A(1)>1$, and suppose that \eqref{rank} fails. With
$\psi$ and $\psi_i$ as in Proposition~\ref{excstructure}, put $\psi_N=1$ and let $\mathcal V_i$
be the closed interval with end points $\psi_i$ and $\psi_{i+1}$. Assume
$0\le\psi_i\le1$ for all $i$, and let $\delta=D_1-1\in(0,1)$, so that $\sum_ia_ib_i^{\delta}=1$. If
\begin{equation}\label{balancedeq}
 \sum_{i:\,y\in \mathcal V_i}b_i^{\delta}=1\qquad\text{for almost every }y\in(0,1),
\end{equation}
then the Hausdorff and both Minkowski dimensions of $G$ equal $D_1$. For equal widths
\eqref{balancedeq} says that almost every level $y$ lies in the same number of intervals
$\mathcal V_i$. This number is then $A(1)$, and $D_1=1+\log A(1)/\log N$.
\end{proposition}
\begin{proof}
The maps $\tau_i(x)=\psi_i+s_ix$ send $[0,1]$ onto $\mathcal V_i\subset[0,1]$, so $\psi$ takes
values in $[0,1]$. Let $\mu$ be the self-affine measure on the graph of $\psi$ with weights
$p_i=a_ib_i^{\delta}$ and let $\nu$ be its image under $(t,x)\mapsto x$. Then
$\nu=\sum_ip_i\,\tau_{i*}\nu$, and this equation has only one solution among probability
measures on $[0,1]$. Lebesgue measure $\mathrm{Leb}$ on $[0,1]$ is a solution:
$\tau_{i*}\mathrm{Leb}$ has density $a_i^{-1}\mathbf 1_{\mathcal V_i}$, and
$\sum_ip_ia_i^{-1}\mathbf 1_{\mathcal V_i}=\sum_ib_i^{\delta}\mathbf 1_{\mathcal V_i}=1$ almost everywhere by
\eqref{balancedeq}. Hence $\nu=\mathrm{Leb}$. For a word $w$ the part of $\mu$ on the cylinder
of $w$ is $p_w$ times the image of $\mu$ under a diagonal affine map with second coordinate
$\tau_w$. Writing $\mathcal V_w=\tau_w([0,1])$, its image under $(t,x)\mapsto x$ is therefore
$p_wa_w^{-1}\mathrm{Leb}|_{\mathcal V_w}=b_w^{\delta}\mathrm{Leb}|_{\mathcal V_w}$.

Let $0<u<1$ and take the words $w$ with $b_w<u\le b_{w^-}$, where $w^-$ is $w$ without
its last symbol. Their time intervals tile $[0,1]$ and have length at least $u\min_ib_i$,
so a ball $\mathcal B$ of radius $u$ in the maximum norm meets at most $C=2+2/\min_ib_i$ of them.
On each of them $\mathcal B$ has mass at most $b_w^{\delta}\cdot2u\le2u^{1+\delta}$. Thus
$\mu(\mathcal B)\le2Cu^{1+\delta}$, and the mass distribution principle gives
$\dim_{\rm H}G\ge1+\delta=D_1$. The upper bound is that of Theorem~\ref{main}. For equal widths
$b_i^{\delta}=N^{-\delta}$ for all $i$, and the integral over $[0,1]$ of the number of intervals
containing $y$ is $\sum_i|\mathcal V_i|=A(1)$.
\end{proof}

\begin{example}\label{balancedex}
Take six equal widths, the factors $(1/4,\,3/4,\,-1/2,\,-1/2,\,2/3,\,1/3)$ and the interpolation
values $0,1/4,1,1/2,0,2/3,1$, their partial sums. Then $\Delta_i=s_i-1/6$, so \eqref{rank}
fails with $\vartheta=1$ and $\psi=f$. The polygonal path through the interpolation values traverses $[0,1]$ three times. Thus almost every level belongs to three intervals $\mathcal V_i$, and
$A(1)=3$. Hence
\[
 \dim_{\rm H}G=\dim_{\B}G=1+\frac{\log 3}{\log 6}\approx1.6131 .
\]
Subdividing a monotone passage at any interior level gives a
family of examples with factors of different moduli. For Okamoto's functions~\cite{okamoto} with
$1/2<a<1$, levels in $[1-a,a]$ belong to three intervals and all other levels in $[0,1]$ to one, so
\eqref{balancedeq} fails.
\end{example}

\subsection{Zero factors and dimension gaps}

\begin{proposition}[Zero branches]\label{zerobranches}
Suppose that $s_i=0$ for at least one $i$, and let $K_J\subset G$ be the attractor of the active branch maps, with $K_J=\varnothing$ when $J=\varnothing$. Then
\[
 \dim_{\rm H}G=\max\{1,\dim_{\rm H}K_J\}\qquad(\dim_{\rm H}\varnothing=0).
\]
\end{proposition}
\begin{proof}
If $s_i=0$, the part of $G$ over $I_i$ is a line segment, and the branch maps are affine.
A point of $G$ that has a coding containing a symbol corresponding to a zero branch therefore lies in the image
of such a segment under a finite composition of active branch maps. There are countably many such segments, at least one of which is nondegenerate. The points that admit an all-active coding form $K_J$. Their union is $G$, so countable stability of Hausdorff dimension gives the formula.
\end{proof}

\paragraph{Example 1. A gap with all branches active.}\hypertarget{gap-active}{}
Take five branches of equal width. Choose interpolation values $(0,1/3,0,1/3,2/3,1)$ and
vertical factors $(1/3,-1/3,1/3,1/3,1/3)$. Then $\Delta_i=s_i-1/5$ and $W=0$.
Each level-$m$ column has width $5^{-m}$ and vertical range of length
$3^{-m}$, since the vertical maps $y\mapsto y_i+s_i y$ cover $[0,1]$.
Counting squares of side $5^{-m}$ gives $\dim_{\B}G=2-\log_5 3$.
The row counts play the role of fibre counts in the carpets of Bedford and McMullen~\cite{bedthesis,mcm}.
Use the two states $G$ and $R(G)$, where $R(t,y)=(t,1-y)$. In the $3$-by-$5$ grid each state contains one orientation-preserving image per column. The negative branch switches states, so the system is irreducible. For each ternary row digit, the matrix below counts branches from the state indexing the row to the state indexing the column.
\[
 M_0=\begin{pmatrix}2&1\\0&1\end{pmatrix},\qquad
 M_1=\begin{pmatrix}1&0\\0&1\end{pmatrix},\qquad
 M_2=\begin{pmatrix}1&0\\1&2\end{pmatrix}.
\]
For a row word $w$, let $c_+(w),c_-(w)$ be the two row sums of $M_{w_1}\cdots M_{w_k}$.
Both columns of $M_j$ have the same sum, equal to $2,1,2$ for $j=0,1,2$, respectively. Hence $c_+(w)+c_-(w)=2\cdot 2^{u(w)}$, where $u(w)$ counts digits different from $1$.
Reflection interchanges $0$ and $2$, so the two states give the same sum of $\alpha$th powers over all words of length $k$, with $\alpha=\log_5 3$.
Since $0<\alpha<1$, subadditivity and monotonicity of $x^\alpha$ give
\[
 \tfrac12\sum_{|w|=k}(2\cdot 2^{u(w)})^\alpha
 \ \le\ \sum_{|w|=k}c_+(w)^\alpha
 \ \le\ \sum_{|w|=k}(2\cdot 2^{u(w)})^\alpha .
\]
The sum on the right is $2^\alpha(2^\alpha+1+2^\alpha)^k$.
Apply the row formula of Kenyon--Peres and Takahashi~\cite{kp,tak} in the form of Allaart and Kawamura~\cite[Equation~(7)]{ak} with $X_1=G$ and $X_2=R(G)$. Starting from $G$, each column word determines one path. Its cylinder meets the interior of its row rectangle because it is a copy of a continuous graph joining the two horizontal edges. Thus the number of open rectangles in row $w$ meeting $G$ is exactly $c_+(w)$. The formula gives
\[
 \HH G=\lim_{k\to\infty}\frac1k\log_3\sum_{|w|=k}c_+(w)^\alpha
 =\log_3(2^\alpha+1+2^\alpha)<2-\alpha=\dim_{\B}G .
\]
The strict inequality follows from strict concavity at $2,1,2$. For dimension results on graphs of continuous
self-affine functions, see K\^ono and Urba\'nski~\cite{kono1,urb}.

\paragraph{Example 2. A gap with zero branches.}\hypertarget{gap-zero}{}
Take four intervals of width $1/4$. The interpolation values are $(0,0,0,1,0)$ and the vertical factors are $(4/5,7/10,0,0)$.
Points whose coding contains a zero branch lie in a countable union of affine images of segments.
Points with an all-active coding have vertical coordinate zero, since both active vertical maps fix zero.
Hence $\HH G=1$. The oscillations over the active level-$m$ columns sum
to $(3/2)^m\operatorname{osc}f$. Counting squares of side $4^{-m}$
and applying the upper bound in Theorem~\ref{main}(a) gives
$\dim_{\B}G=1+\log_4(3/2)$.

\begin{remark}[Scope of the exceptional set]\label{excremarks}
\textup{(i)} \emph{Okamoto's functions.} For $N=3$, equal widths, interpolation values $0,a,1-a,1$ and
factors $(a,1-2a,a)$ the interpolant is the function introduced by Okamoto~\cite{okamoto}. Here
$\Delta_i=s_i-1/3$. For $a\ne1/3$ the residual rank is one and \eqref{rank} fails with $\vartheta=1$. At $a=1/3$ the interpolant is the identity and \eqref{rank} holds at rank zero. For $1/2<a<1$ Prokaj and
B\'ar\'any~\cite[Theorem~1.1]{pb} prove that the Hausdorff, box and Assouad dimensions of the graph equal
$1+\log_3(4a-1)=D_1$ for all $a$ outside a set of Hausdorff dimension zero.

\textup{(ii)} \emph{Exact overlaps.} In the \hyperlink{gap-active}{five-branch example} above,
$\psi=f$ and the maps $x\mapsto\psi_0+s_0x$ and $x\mapsto\psi_2+s_2x$ coincide. In that
example the Hausdorff dimension is strictly smaller than $D_1$, and the weighted covering
numbers of Proposition~\ref{balanced} are $9/5,3/5,3/5$ on the interiors of the three thirds of $[0,1]$. Exact
overlaps alone do not force a drop: in Example~\ref{balancedex} the maps $x\mapsto x/4$ and
$x\mapsto 2x/3$ commute, and the three dimensions agree.

\textup{(iii)} \emph{What remains open.} For all-active vectors in $\Xi$ with $A(1)>1$ the
general problem remains open for $r=1$ outside the sufficient conditions above and for $r\ge2$
with $A(r-1)>1$. In the first case the question is the dimension of a planar diagonal
self-affine graph, which is related to the self-similar measures of $\mathcal F_\psi$. A formula
valid for all parameters would in particular settle the dimension of every Okamoto function.
In the second case $D_{r-1}<D_r$. We expect that $\dim_{\rm H}G=D_r$ when $\mathcal F_\psi$ is
exponentially separated. The proof of Theorem~\ref{hom}\textup{(1)} uses the hypothesis on the
shear span in one place, the separation of projected translations by the slope
(Lemma~\ref{sep}), and there only for one row of $P$. For block systems obtained from
exceptional data the shear differences span $\ker\vartheta$ and the stationary slope measure is supported on an
affine hyperplane. Any row of $P$ that is nonzero on $\ker\vartheta$ can be used in the proof of Lemma~\ref{sep}. Since the rows are orthonormal, the only remaining case is a single row proportional to $\vartheta$. At the invariant slope, the projected measure is self-similar for a homogeneous system on the line. An extension of Theorem~\ref{hom} to this situation
needs a hypothesis on that system, and we do not pursue it here.

\textup{(iv)} \emph{Zero branches.} If the residual increments on the active branches do not span $E$,
condition \eqref{rank} fails for all values of the nonzero factors, as in the \hyperlink{gap-zero}{four-branch example} above.
\end{remark}

\section{Two independent counterexamples to Massopust's formula}\label{masssec}
For two diagonal vertical coordinates with coefficients $e_{i,1},e_{i,2}$ and noncoplanar interpolation points, part (a) of \cite[pp.~8--9]{mass} asserts that the Minkowski dimension is the solution $d_{\mathrm M}$ of $\sum_i |e_{i,1}e_{i,2}|b_i^{d_{\mathrm M}-2}=1$ when the sum of the products exceeds one. Part (c) concerns the case in which the sum of the products is at most one while both coordinate sums exceed one. As printed, (3.3)--(3.4) of~\cite{mass} then give the larger of the two exponents $d_k'$ defined by $\sum_i|e_{i,k}|b_i^{d_k'-n}=1$, $k=1,2$, where $n=2$ is the number of vertical coordinates. The proof for $n=2$ (see (3.14) and pp.~17--18 there) works instead with the single-coordinate exponents $d_k$ given by $\sum_i|e_{i,k}|b_i^{d_k-1}=1$ and arrives at $\max\{d_1,d_2\}$. We read part (c) in this sense. Our example for part (c) also contradicts the printed exponent. Ma\l{}ysz~\cite[p.~418]{malysz2001} identified part~(c) of the main theorem of~\cite{mass} as giving an incorrect Minkowski dimension for graphs of vector-valued FIFs and referred to a manuscript then in preparation for details.

Start a lattice walk at $v_0=(0,0)$, end at $v_N=(3,3)$, and use steps in $\{-1,1\}^2$. Set $P_j=(j/N,v_j/3)$ and
\[
 \Phi_i z=\tfrac13(v_{i-1}+\operatorname{diag}(v_i-v_{i-1})z),\qquad
 \mathcal M_i(t,z)=((i-1+t)/N,\Phi_i z).
\]
These maps join consecutive interpolation points: $\mathcal M_i$ sends $P_0$ to $P_{i-1}$ and $P_N$ to $P_i$. Below, $G$ is the attractor of $\mathcal M_1,\ldots,\mathcal M_N$. It is the graph of the continuous interpolant of the interpolation points $P_j$ with equal widths $b_i=1/N$ and diagonal vertical parts, whose entries $e_{i,k}\in\{-1/3,1/3\}$ are the coordinates of $(v_i-v_{i-1})/3$. Its vertical range is the attractor of the $\Phi_i$. All vertical moduli are $1/3$. We write $d_1,d_2$ for the exponents of part (c) in the reading above, $d_1',d_2'$ for the printed ones, and $d_{\mathrm M}$ for the exponent of part (a).

\begin{proposition}[Counterexample to part (c)]\label{prop:c}
Let $N=9$ and take the classical Peano parametrization of the square, in the orientation of Allaart and Kawamura~\cite[Section~4.1, Figure~2]{ak}. Its walk has vertices
\[
 (0,0),(1,1),(0,2),(1,3),(2,2),(1,1),(2,0),(3,1),(2,2),(3,3).
\]
Then the interpolation points $P_j$ are not coplanar, $\sum_i|e_{i,1}e_{i,2}|=1$, $\sum_i|e_{i,k}|=3$ for $k=1,2$, and
\[
 \HH G=\underline{\dim}_{\B}G=\overline{\dim}_{\B}G=2 ,
\]
whereas $\max\{d_1,d_2\}=3/2$ and $\max\{d_1',d_2'\}=5/2$.
\end{proposition}
\begin{proof}
In the indexing of the maps $\Phi_1,\ldots,\Phi_9$, the first coordinate reverses orientation on branches $2,5,8$ and the second on $4,5,6$. With the prescribed endpoints these are precisely the two coordinate recursions of the Peano curve. The nine maps visit each cell of the $3\times3$ square grid exactly once, so their images tile the square. Hence projection gives $\HH G\ge2$. Each of $9^k$ pieces fits in a cube of side $3^{-k}$, hence $\overline{\dim}_{\B}G\le2$. The determinant of the three differences of interpolation points from $P_0$ with indices $1,2,4$ is $4/81$. The product sum equals one and each coordinate sum is three, so part (c) applies. Since $b_i=1/9$ and $|e_{i,k}|=1/3$, the equation for $d_k$ reads $3\cdot9^{1-d_k}=1$. It gives $3/2$ instead of the actual value $2$. The equation for the printed exponent reads $3\cdot9^{2-d_k'}=1$. It gives $5/2$, which is not the actual value either.
\end{proof}

\begin{proposition}[Counterexample to part (a)]\label{prop:a}
Let $N=729=3^6$, let $v_j=(j\bmod2,j\bmod2)$ for $0\le j\le723$, and let $v_{724},\ldots,v_{729}$ be
\[
 (0,2),(1,3),(2,2),(1,1),(2,2),(3,3).
\]
Then the interpolation points $P_j$ are not coplanar, $\sum_i|e_{i,1}e_{i,2}|=81$, $d_{\mathrm M}=8/3$, and
\[
 \overline{\dim}_{\B}G\le1+\tfrac56\log_3 8<8/3.
\]
\end{proposition}
\begin{proof}
The determinant using indices $1,2,724$ is $-4/6561$. The product sum is $729/9=81$, and $81\cdot729^{2-d_{\mathrm M}}=1$ gives $d_{\mathrm M}=8/3$. Only two distinct oriented vertical maps occur in the first $723$ branches and at most six more occur in the tail. In each horizontal prefix of length $k$ refine vertically by $5k$ further symbols. The $8^{5k}$ resulting squares have side $3^{-6k}$, also the prefix time width. Thus $729^k8^{5k}$ cubes of side $3^{-6k}$ cover $G$, which gives the bound $1+\tfrac56\log_38$ along the scales $3^{-6k}$. The geometric sequence suffices by monotonicity. Keeping orientations in this count is essential. The strict inequality holds because $\log_38<2$.
\end{proof}
Both examples have full affine rank. The lower-bound argument in \cite[equations (3.10)--(3.13)]{mass} uses the product of the covering counts for the two coordinate projections. A diagonal segment shows why multiplying the two counts can overestimate the number of squares needed in the plane. Each coordinate projection meets on the order of $\varepsilon^{-1}$ grid intervals, while the segment itself meets only on the order of $\varepsilon^{-1}$ grid squares.

\paragraph{Part (b) of Massopust's theorem.} When every $\sum_i|e_{i,j}|\le1$ the Minkowski dimension is one as claimed in \cite{mass}. Iterate the interpolation operator on the endpoint affine interpolant. For each coordinate the variation of its polygonal approximants satisfies $V_{k+1,j}\le C_j+(\sum_i |e_{i,j}|)V_{k,j}$. Thus the graphs have length $O(k+1)$ and converge uniformly at rate $O(\theta^k)$ for any $\max_{i,j}|e_{i,j}|<\theta<1$. With $k$ of order $\log(1/\varepsilon)$ their arc-length covers give $O(\varepsilon^{-1}(1+\log(1/\varepsilon)))$ cubes for the limit. Projection gives the matching dimension.

\begin{remark}[Size of the discrepancy in part (c)]\label{prop:disc}
For the positive three-branch example in Table~\ref{tab:threebranch}, take
$1/3<s\le1/\sqrt3$ and put $q=\log3/\log(1/s)$. The product sum is $3s^2\le1$ and both coordinate sums equal $3s>1$. The two single-coordinate equations in the proof of part (c) give $d_1=d_2=2+\log_3s=2-1/q$. The printed exponent is $3+\log_3s$. Theorem~\ref{main} gives graph dimension $q$, so the discrepancy from $d_1$ is
\[
 q-2+1/q=(q-1)^2/q .
\]
It increases on $1<q\le2$ and reaches $1/2$ at $s=1/\sqrt3$. This is the maximum within the given family over the stated parameter interval.
\end{remark}

\begin{remark}[Failure of the projected affine dimension bound for diagonal maps]\label{rem:diag-bound}
Theorem~\ref{main} treats scalar vertical actions on the residual space. For dimension results on diagonal self-affine systems see Rapaport~\cite{rap}.

The bound by $g$ in Theorem~\ref{main}(b) can fail for positive diagonal vertical parts. Take two equal-width branches with $f(0)=f(1)=(0,0)$, $f(1/2)=(1,1)$, and vertical part $\operatorname{diag}(3/5,4/5)$ on both branches. For $\Pi(t,z_1,z_2)=z_2-z_1$, all interpolation points project to zero, yet $f(1/4)=(11/10,13/10)$, so $\Pi(G)$ contains $[0,1/5]$.
\end{remark}

\appendix
\section{Self-affine measures on the graph and their tilted projections}\label{measureappendix}
We prove the two measure assertions of Theorem~\ref{hom}. Throughout this appendix logarithms have base two.
\subsection{Entropy of the graph measure and exact dimensionality}\label{exactsec}
\begin{lemma}[Entropy dimension of the graph measure]\label{entgraph}
The entropy dimension of $\mu_{\boldsymbol\pi}$ exists and
\[
 \lim_{j\to\infty}\frac{H_j(\mu_{\boldsymbol\pi})}{j}=\hw/X+(1-\ell_S/X)\min(r,Q).
\]
\end{lemma}
\begin{proof}
Take $P=I_r$, let $W$ be the $\boldsymbol\pi$-distributed length-$k$ prefix, and put $j=\lceil kX\rceil$, so $2^{-j}\le\rho^k<2^{1-j}$. A time cell meets at most three prefix intervals and each interval meets at most three time cells, including endpoints. Hence $H(W\mid\mathcal D_j(T,Z))\le\log3$ and $H(\mathcal D_j(T)\mid W,\mathcal D_j(Z))\le\log3$. The chain rule applied to $(W,\mathcal D_j(T),\mathcal D_j(Z))$ gives
\[
 |H_j(\mu_{\boldsymbol\pi})-k\hw-H(\mathcal D_j(Z)\mid W)|\le\log3.
\]
Given $W=w$, \eqref{cyl} with slope zero rescales the vertical mesh to $2^{-j}S^{-k}$. This differs by a factor of at most four from $2^{-\lfloor k(X-\ell_S)\rfloor}$. Lemma~\ref{grid}(2) gives
\[
 \left|H_j(\mu_{\boldsymbol\pi})-k\hw-\E_wH_{\lfloor k(X-\ell_S)\rfloor}(\nu_{R_w(0)})\right|
 \le\log3+C_{\rm grid}(r,4).
\]
Here $R_w(0)\in B$. Divide by $j$ and apply Theorem~\ref{hom}(1) with $P=I_r$, noting that $k/j\to1/X$ and $\lfloor k(X-\ell_S)\rfloor/j\to1-\ell_S/X$. The resulting limit along the scales $j$ is $\hw/X+(1-\ell_S/X)\min(r,Q)$. The limit extends to all integer scales because consecutive values of $\lceil kX\rceil$ differ by at most $\lceil X\rceil$ and the entropy increases by at most $r+1$ per step.
\end{proof}

\begin{proof}[Proof of Theorem~\ref{hom}(2)]
Conjugating by $(t,z)\mapsto(t,az)$ for sufficiently small $a>0$ makes all branches contractions in the Euclidean norm. Feng's theorem for the coding image of an ergodic Bernoulli measure~\cite[Theorem~1.2(ii)]{feng} gives exact dimensionality, also for overlapping images. Bi-Lipschitz conjugacy preserves local dimension. For self-affine systems with a totally dominated splitting see also B\'ar\'any and K\"aenm\"aki~\cite[Theorem~2.6]{bk}.

A fixed translation and homothety place the compact support of $\mu_{\boldsymbol\pi}$ inside the unit cube. This similarity preserves local dimensions and changes $H_j$ by $O(1)$ uniformly in $j$, by Lemma~\ref{grid}(2). For an exact-dimensional probability measure on the cube, the common local dimension equals the limit of its normalized dyadic entropies~\cite[Section~6.4, equation~(73)]{hoch} (see also \cite[Theorem~4.4]{young82}). Lemma~\ref{entgraph} computes this limit for $\mu_{\boldsymbol\pi}$ and gives the claimed value.
\end{proof}

\subsection{Local dimensions of projected measures}
\begin{proof}[Proof of Theorem~\ref{hom}(3)]
Fix $P,a$ and write $Q_p=\min(p,Q)$. We first bound the lower local dimension using conditional entropies along the coding tree, then bound the upper local dimension using cylinder masses and the ambient dimension. Let $\mathbb P$ be the Bernoulli product measure on symbol sequences and
$F(\omega)=PZ(\omega)+aT(\omega)$ the projected coding map.
Choose an invariant slope ball containing $a$ and let $B_Y$ bound all the corresponding
projected supports. A cylinder with prefix $w$ of length $k$ has projected diameter at most
$2B_Y S^k$. The corresponding conditional projected measure is a translate of $S^k\nu_{R_w a}$ by \eqref{cyl}.

Group the symbols into blocks of length $q$ and give the resulting symbolic tree the metric
$(S^q)^j$ for paths with precisely $j$ common initial blocks. The Lipschitz constant of $F$
is at most $2B_Y$, independently of $q$. For a probability measure $\theta$ put
$\mathcal H_s(\theta)=-\int\log_2\theta(B(y,s))\,d\theta(y)$.
For every prefix of $j$ blocks, the $(S^q)^{j+1}$-entropy of its conditional image
is exactly $\mathcal H_{S^q}(\nu_{R_w a})$. This differs from
$H_{\lfloor q\ell_S\rfloor}(\nu_{R_w a})$ by a bounded dimensional constant
\cite[Lemmas 3.3 and 3.5]{hsaverages}. For every $\varepsilon>0$, part (1) therefore
bounds these conditional entropies below by $(Q_p-\varepsilon)q\ell_S-O(1)$ for all large $q$,
uniformly in $j,w$. The local entropy averages theorem on a tree now gives
$\underline{\dim}_{\rm loc}(\nu_a,y)\ge Q_p-\varepsilon-O(1/q)$ for $\nu_a$-almost every $y$
\cite[Theorems 4.4 and 5.4(4), Proposition 5.2]{hsaverages}; see also
\cite[Lemma 4.8]{algomshmerkin}. Here one block changes the logarithmic scale by $q\ell_S=q\log_2(1/S)$, so a bounded entropy error per block contributes $O(1/q)$ to the dimension bound. The error constant depends on the target dimension,
$S$ and the common Lipschitz bound. Taking countably many block lengths and
then $\varepsilon\downarrow0$ gives the lower bound $Q_p$.

For $\mathbb P$-almost every $\omega$, the strong law of large numbers gives
$-n^{-1}\log_2\mathbb P[\omega|n]\to H(\boldsymbol\pi)$.
The ball of radius $(1+2B_Y)S^n$ about $F(\omega)$ contains the entire projected
$n$-cylinder, and hence has mass at least $\mathbb P[\omega|n]$.
Consecutive radii have fixed ratio $S$, so the upper local dimension is at most $Q$
at $\nu_a$-almost every point. It is also at most $p$: for $\gamma>0$, the total mass
of dyadic cubes of side $2^{-n}$ having mass below $2^{-n(p+\gamma)}$
is $O(2^{-n\gamma})$ by bounded support. Borel--Cantelli and containment of each
such scale cube in a ball of radius $(1+\sqrt p)2^{-n}$ about any of its points give
the bound $p+\gamma$ almost everywhere. Let $\gamma\downarrow0$.
The two bounds give the asserted limit for this $P,a$.
\end{proof}

\section*{Acknowledgements}
From September through November 2026, I used
\href{https://chatgpt.com}{ChatGPT (Astra, OpenAI)} and
\href{https://claude.ai}{Claude (Fable 5.2, Anthropic)} to explore proof strategies,
test conjectures and search for counterexamples. I used both tools in developing
parts of the direct proof of the lower bound for the graph dimension, checking intermediate
estimates and formulating Theorem~\ref{main} for all linear images of the graph.
They also assisted with the Lean formalization, English editing, \LaTeX{}, and reference formatting.
I take full responsibility for the mathematical claims and their proofs, the formalization, the references, and the final manuscript.

\section*{Competing interests}
The author declares no competing interests.

\end{document}